\documentclass[12pt,a4paper]{article}
\usepackage{multirow}
\usepackage[dvips]{graphicx}
\usepackage{latexsym,amsmath,amsfonts,amscd}
\usepackage{epsfig}
\usepackage{indentfirst}
\usepackage{verbatim}
\usepackage{color}
\usepackage{colordvi}
\usepackage{tikz}
\usepackage{array}

\begin{document}

\baselineskip=2pc

\begin{center}
  {\Large \bf A  Hermite WENO scheme with artificial linear weights for hyperbolic conservation laws\footnote{The research is partly supported by Science Challenge Project, No. TZ2016002  and NSAF grant U1630247.} }
\end{center}

\centerline{Zhuang Zhao\footnote{School of Mathematical Sciences, Xiamen University,
Xiamen, Fujian 361005, P.R. China. E-mail: zzhao@stu.xmu.edu.cn.} and Jianxian
Qiu\footnote{School of Mathematical Sciences and Fujian Provincial
Key Laboratory of Mathematical Modeling and High-Performance
Scientific Computing, Xiamen University,
Xiamen, Fujian 361005, P.R. China. E-mail: jxqiu@xmu.edu.cn.}}

\vspace{.1in}

\centerline{\bf Abstract}

\bigskip

In this paper, a fifth-order Hermite weighted essentially non-oscillatory (HWENO) scheme with artificial linear weights is proposed for one and two dimensional  hyperbolic conservation laws, where the zeroth-order and the first-order moments are used in the spatial reconstruction.  We  construct the HWENO methodology using a nonlinear convex combination of a high degree polynomial with several low degree polynomials, and the associated linear weights can be any artificial positive numbers with only requirement that their summation equals one. The one advantage of the  HWENO scheme is its simplicity and easy extension to multi-dimension in engineering applications for we can use any artificial linear weights which are independent on geometry of mesh.  The another advantage is its higher order numerical  accuracy using less candidate stencils for two dimensional problems. In addition, the  HWENO scheme still keeps the compactness as only immediate neighbor information is needed in the reconstruction and has high efficiency for directly using linear  approximation in the smooth regions. In order to avoid nonphysical oscillations nearby strong shocks or contact discontinuities, we adopt the thought of limiter for discontinuous Galerkin method to control the spurious oscillations. Some benchmark numerical tests are performed to demonstrate the capability of the proposed scheme.

\vfill {\bf Key Words:} Hermite WENO scheme, hyperbolic conservation laws, unequal size spatial stencil, hybrid, discontinuous Galerkin method

{\bf AMS(MOS) subject classification:} 65M60, 35L65

\pagenumbering{arabic}
\newpage

\section{Introduction}
\label{sec1}
\setcounter{equation}{0}
\setcounter{figure}{0}
\setcounter{table}{0}

In this paper, we develop a fifth order  Hermite weighted essentially non-oscillatory (HWENO) scheme with artificial linear weights for one and two dimensional nonlinear hyperbolic conservation laws. The idea of  HWENO scheme is similar to that of weighted essentially non-oscillatory (WENO) scheme which have been widely applied for computational dynamics fluids.   In 1994, the first WENO scheme was proposed by Liu, Osher and Chan \cite{loc}  mainly in terms of ENO scheme \cite{ho,h1,heoc}, in which they combined all candidate stencils by a nonlinear convex manner to obtain higher order accuracy in smooth regions, then, in 1996, Jiang and Shu \cite{js} constructed the third and fifth-order finite difference WENO schemes in multi-space dimension, where they gave a general definition  for smoothness indicators and nonlinear weights. Since then, WENO schemes have been further developed in \cite{hs,lpr,SHSN,ccd,ZQd}.
However, if we design a higher order accuracy WENO scheme, we need to enlarge the  stencil. In order to keep the compactness of the scheme, Qiu and Shu \cite{QSHw1,QSHw} gave a new option by evolving both with the solution and its derivative, which were termed as Hermite WENO (HWENO) schemes.

HWENO schemes would have higher order accuracy than  WENO schemes with the same reconstruction stencils.  As the solutions of nonlinear hyperbolic conservation laws often contain discontinuities, its derivatives or first order moments would be relatively large nearby discontinuities. Hence, the HWENO schemes presented in \cite{QSHw1,QSHw,ZQHW,TQ,LLQ,ZA,TLQ,CZQ} used  different  stencils for discretization in  the space for the original  and  derivative equations, respectively.  In one sense, these HWENO schemes can be seen as an extension by DG methods, and Dumbser et al. \cite{DBTM} gave a general and unified framework to define the numerical scheme extended by DG method, termed as $P_NP_M$ method. But the derivatives or the first order moments  were still used straightforwardly nearby the discontinuities, which would be less robust for problems with strong shocks. Such as the first HWENO schemes \cite{QSHw1,QSHw} failed to simulate the double Mach and the step forward problems, then, Zhu and Qiu \cite{ZQHW} solved this problem by using a new procedure to reconstruct the derivative terms, while Cai et al. \cite{CZQ} employed additional positivity-preserving manner. Overall, only using different stencils to discretize the space is not enough to overcome the effect of the derivatives or the first order moments near the discontinuities.  Hence, we took the thought of limiter for discontinuous Galerkin (DG) method \cite{DG2}  to modify the first order moments  nearby the discontinuities of the solution in \cite{ZCQHH}, meanwhile, we also noticed that many hybrid WENO schemes \cite{SPir,Hdp,cd1,cd2,Glj,ZZCQ} employed linear schemes directly in the smooth regions, while still used WENO schemes in the discontinuous regions, which can increase the efficiency obviously, therefore, in \cite{ZCQHH}, we directly used high order linear approximation in the smooth regions, while modified the first order moments on the troubled-cells and employed HWENO reconstruction on the interface. The hybrid HWENO scheme \cite{ZCQHH} had high efficiency and resolution with non-physical oscillations, but it still had a drawback of that the linear weights were depended on geometry of the mesh and point where the reconstruction was performed,  and they were not easy to be computed, especially for multi-dimensional problems with unstructured meshes. For example, in \cite{ZCQHH} we needed to compute the linear weights at twelve points in one cell by a least square methodology with eight small stencils for two dimensional problems, in which the numerical accuracy was only the fourth order.  Moreover, if we solve the problems for unstructured meshes, the linear weights would be more difficult to calculate, and the negative weights may appear or there is nonexistent of the linear weights for some cases.  In order to overcome the drawback,  Zhu and Qiu \cite{WZQ} presented a new simple WENO scheme in the finite difference framework, which had a convex combination of a fourth degree polynomial and other two linear polynomials by using any artificial positive linear weights (the sum equals one). Then the method was extended to finite volume methods both in structured and unstructured meshes \cite{bgs,vozq, DBSR, tezq,trzq,bgfb}.

In this paper, following the idea of the new type WENO schemes \cite{WZQ, vozq, DBSR, tezq, trzq}, hybrid WENO \cite{SPir,Hdp,cd1,cd2,Glj,ZZCQ} and hybrid HWENO \cite{ZCQHH},  we develop the new hybrid HWENO scheme in which we use a nonlinear convex combination of a high degree polynomial with several low degree polynomials and the linear weights can be any artificial positive numbers with the only constraint that their sum is one.
The procedures of the new hybrid HWENO scheme are: firstly, we modify the first order moments using the new HWENO limiter methodology in the troubled-cells, which are identified by the KXRCF troubled-cell indicator \cite{LJJN}. Then, for the space discretization, if the cell is identified as a troubled-cell, we would use the new HWENO reconstruction at the points on the interface; otherwise we employ linear approximation at the interface points straightforwardly. And we directly use high order linear approximation at the internal points for all cells. Finally, the third order TVD Runge-Kutta method \cite{so1} is applied for the time discretization. Particularly, only the new HWENO reconstructions need to be performed on local characteristic directions for systems. In addition, the new hybrid HWENO scheme inherits the advantages of \cite{ZCQHH}, such as non-physical oscillations for using the idea of limiter for discontinuous Galerkin (DG) method, high efficiency for employing linear approximation straightforwardly in the smooth regions, and compactness as only immediate neighbor information is needed, meanwhile, it gets less numerical errors on the same meshes and has higher order numerical accuracy for two dimensional problems.

The organization of the paper is as follows: in Section 2, we introduce the detailed implementation of  the new hybrid HWENO scheme in the one and two  dimensional cases. In Section 3, some benchmark numerical are performed to illustrate the numerical accuracy, efficiency, resolution and robustness of proposed scheme. Concluding remarks are given in Section 4.

\section{Description of  Hermite WENO scheme with artificial linear weights}
\label{sec2}
\setcounter{equation}{0}
\setcounter{figure}{0}
\setcounter{table}{0}
In this section, we present the construction procedures of the hybrid HWENO scheme with artificial linear weights for  one and two dimensional hyperbolic conservation laws, which is the fifth order accuracy both in the one and two dimensional cases.

\subsection{One dimensional case}
At first, we consider one dimensional scalar hyperbolic conservation laws
  \begin{equation}
\label{EQ} \left\{
\begin{array}
{ll}
u_t+ f(u)_x=0, \\
u(x,0)=u_0(x). \\
\end{array}
\right.
\end{equation}
The computing domain is divided by uniform meshes $I_{i} =[x_{i-1/2},x_{i+1/2}]$ for simplicity,  the mesh center $x_i=\frac {x_{i-1/2}+x_{i+1/2}} 2$ with the mesh size  $\Delta x=x_{i+1/2}- x_{i-1/2}$.

As the variables of our designed  HWENO scheme are the zeroth and first order moments, we multiply the governing equation (\ref{EQ}) by $\frac{1}{\Delta x}$ and $\frac{x-x_i}{{(\Delta x)}^2}$, respectively, and integrate  them over $I_i$, then, employ the numerical flux to approximate the values of the flux at the interface. Finally, the semi-discrete finite volume HWENO scheme is
\begin{equation}
\label{odeH}
  \left \{
  \begin{aligned}
   \frac{d \overline u_i(t)}{dt} &=- \frac 1 {\Delta x} \left ( \hat f_{i+1/2}-\hat f_{i-1/2}\right ),\\
    \frac{d \overline v_i(t)}{dt}&=- \frac 1 {2\Delta x} \left ( \hat f_{i-1/2} +\hat f_{i+1/2}\right ) +\frac 1 {\Delta x} F_i(u).
  \end{aligned}
  \right.
\end{equation}
The initial conditions are $\overline u_i(0) = \frac 1 {\Delta x} \int_{I_i} u_0(x) dx$ and $\overline v_i(0) = \frac 1 {\Delta x} \int_{I_i} u_0(x) \frac{x-x_i} { \Delta x} dx$. $\overline u_i(t)$ is the zeroth order moment in $I_i$ as $\frac 1 {\Delta x} \int_{I_i} u(x,t) dx$ and $\overline v_i(t)$ is the first order moment in $I_i$ as $\frac 1 {\Delta x} \int_{I_i} u(x,t) \frac{x-x_i} { \Delta x}dx $. $\hat f_{i+1/2}$ is the numerical flux to approximate the value of  the flux $ f(u) $ at the interface point $x_{i+1/2}$, which is defined by the Lax-Friedrichs numerical flux method, and the explicit expression is
\begin{equation*}
\label{Nflux} \hat f_{i+1/2}= \frac 1 2 \left (f(u^-_{i+1/2})+f(u^+_{i+1/2}) \right) -\frac \alpha 2\left( u^+_{i+1/2}-u^-_{i+1/2} \right),
\end{equation*}
in which $\alpha= \max_{u}|f'(u)|$. $F_i(u)$ is the numerical integration for the flux $f(u)$ over $I_i$, and is approximated by a four-point Gauss-Lobatto quadrature formula:
\begin{equation*}
\label{Nfluxinte} F_i(u)=\frac 1 {\Delta x} \int_{I_i} f(u) dx \approx \sum_{l=1}^4 \omega_l f(u(x_l^G,t)),
\end{equation*}
where the weights are $\omega_1=\omega_4=\frac 1 {12}$ and $ \omega_2=\omega_3=\frac 5 {12}$, and the quadrature points on the cell $I_i$ are
\begin{equation*}
x_1^G=x_{i-1/2},\quad x_2^G=x_{i-\sqrt5/10}, \quad x_3^G=x_{i+\sqrt5/10}, \quad x_4^G=x_{i+1/2},
\end{equation*}
in which  $x_{i+a}$ is $x_i +a \Delta x$.

Now, we first present the detailed  procedures of the spatial reconstruction for  HWENO scheme in Steps 1 and 2,  then, we introduce the method of time discretization in Step 3.

\textbf{Step 1.} Identify the troubled-cell and  modify the first order moment in the troubled-cell.

Troubled-cell means that the solution of the equation in the cell may be discontinuous, we first use the KXRCF troubled-cell indicator \cite{LJJN} to identify the troubled-cell, and the procedures were given in the hybrid HWENO scheme \cite{ZCQHH}, then, if the cell $I_i$ is identified as a troubled-cell, we would modify the first order moment $\overline v_i$ by the following procedures.

We use the thought of HWENO limiter \cite{QSHw1} to modify the first order moment, but the modification for the first order moment is based on a convex combination of a fourth degree polynomial with two linear polynomials. Firstly, we give a  large stencil $S_0=\{I_{i-1},I_{i},I_{i+1}\}$ and two small stencils $S_1=\{I_{i-1},I_{i}\}$, $S_2=\{I_{i},I_{i+1}\}$, then, we obtain a quartic polynomial $p_0(x)$ on $S_0$, as
\begin{equation*}
 \frac{1}{\Delta x} \int_{I_{i+j}} p_0(x)dx = \overline u_{i+j},\ j=-1,0,1, \quad \frac{1}{\Delta x} \int_{I_{i+j}} p_0(x)\frac{x-x_{i+j}} {\Delta x} dx = \overline v_{i+j}, \ j=-1,1,
\end{equation*}
and get two linear polynomials $p_1(x),p_2(x)$ on $S_1,S_2$, respectively, satisfying
\begin{equation*}
\begin{split}
\frac{1}{\Delta x} \int_{I_{i+j}} p_1(x)dx &= \overline u_{i+j}, \quad j=-1,0,\\
\frac{1}{\Delta x} \int_{I_{i+j}} p_2(x)dx &= \overline u_{i+j}, \quad j=0,1.
\end{split}
\end{equation*}
We use these three polynomials to reconstruct $\overline v_i=\frac{1}{\Delta x} \int_{I_{i}} u(x)\frac{x-x_{i}} {\Delta x} dx$, and their explicit results  are
\begin{equation*}
\begin{split}
\frac{1}{\Delta x} \int_{I_{i}} p_0(x)\frac{x-x_{i}} {\Delta x}dx &= \frac 5 {76} \overline u_{i+1}-\frac 5 {76} \overline u_{i-1}-\frac {11}{38}\overline v_{i-1}-\frac {11}{38}\overline v_{i+1},\\
\frac{1}{\Delta x} \int_{I_{i}} p_1(x)\frac{x-x_{i}} {\Delta x}dx &=\frac 1 {12} \overline u_i-\frac 1 {12} \overline u_{i-1},\\
\frac{1}{\Delta x} \int_{I_{i}} p_2(x)\frac{x-x_{i}} {\Delta x}dx &= \frac 1 {12}\overline u_{i+1}-\frac 1 {12} \overline u_{i}.\\
\end{split}
\end{equation*}
For simplicity, we define $q_n$ as $\frac{1}{\Delta x} \int_{I_{i}} p_n(x)\frac{x-x_{i}} {\Delta x}dx$ in the next procedures. With the similar idea of the central WENO schemes \cite{lpr,lpr2}  and the new WENO schemes \cite{WZQ,vozq,tezq,trzq}, we rewrite $q_0$ as:
\begin{equation}\label{p_big}
  q_0=\gamma_0 \left( \frac 1 {\gamma_0}q_0 -\frac {\gamma_1} {\gamma_0}q_1- \frac {\gamma_2} {\gamma_0} q_2 \right) + \gamma_1 q_1+ \gamma_2 q_2.
\end{equation}
We can notice that equation (\ref{p_big}) is always satisfied for any choice of $\gamma_0$, $\gamma_1$, $\gamma_2$ with $\gamma_0 \neq 0$. To make the next WENO procedure be stable, the linear weights would be positive with  $\gamma_0+ \gamma_1+\gamma_2=1$, then, we calculate the smoothness indicators $\beta_n$ to measure how smooth the functions $p_n(x)$ in the cell $I_i$, and we use the same definition as in \cite{js},
\begin{equation}
\label{GHYZ}
\beta_n=\sum_{\alpha=1}^r\int_{I_i}{\Delta x}^{2\alpha-1}(\frac{d ^\alpha
p_n(x)}{d x^\alpha})^2dx, \quad n=0,1,2,
\end{equation}
where $r$ is the degree of the polynomials $p_n(x)$, then, the expressions for the smoothness indicators are
\begin{equation*}
  \left \{
  \begin{aligned}
    \beta_0=&(\frac{29}{38} \overline u_{i-1}-\frac{29}{38} \overline u_{i+1} +\frac{60}{19} \overline v_{i-1}+\frac{60}{19} \overline v_{i+1} )^2+( \frac94 \overline u_{i-1}-\frac92 \overline u_{i}+\frac94 \overline u_{i+1} +\frac{15}2\overline v_{i-1}-\frac{15}2\overline v_{i+1}
     )^2+\\
     &\frac{3905}{1444} ( \overline u_{i-1} - \overline u_{i+1}+12 \overline v_{i-1} + 12 \overline v_{i+1}  )^2+\frac{1}{12}( \frac52 \overline u_{i-1}-5\overline u_{i}+\frac52 \overline u_{i+1}+9\overline v_{i-1}-9\overline v_{i+1} )^2+\\
     &\frac{109341}{448}( \overline u_{i-1}-2 \overline u_{i}+ \overline u_{i+1} +\overline v_{i-1}-\overline v_{i+1} )^2,\\
     \beta_1=& (\overline u_i -\overline u_{i-1})^2,\\
     \beta_2=& (\overline u_{i+1} -\overline u_{i})^2.\\
  \end{aligned}
  \right.
\end{equation*}
Later, we use a  new parameter $\tau$ to measure the absolute difference between $\beta_{0}$, $\beta_{1}$ and $\beta_{2}$, which is also can be seen in these new WENO  schemes \cite{WZQ,vozq,tezq,trzq},
\begin{equation}
\label{tao5}
\tau=(\frac{|\beta_{0}-\beta_{1}|+|\beta_{0}-\beta_{2}|}{2})^2,
\end{equation}
and the nonlinear weights are defined as
\begin{equation*}
\label{9}
\omega_n=\frac{\bar\omega_n}{\sum_{\ell=0}^{2}\bar\omega_{\ell}},
\ \mbox{with} \ \bar\omega_{n}=\gamma_{n}(1+\frac{\tau}{\beta_{n}+\varepsilon}), \ n=0,1,2,
\end{equation*}
where $\varepsilon = 10^{-6}$ is to avoid the denominator by zero.
Finally, the  first order moment $\overline v_i$ is modified by
\begin{equation*}
  \overline v_i = \omega_0 \left( \frac 1 {\gamma_0}q_0 - \sum_{n=1}^{2}\frac {\gamma_n} {\gamma_0} q_n\right) + \sum_{n=1}^{2}\omega_n q_n.
\end{equation*}
Noticed that we just replace the linear weights in equation (\ref{p_big}) by the nonlinear weights, and the accuracy of the modification depends on the accuracy of the high degree reconstructed polynomial. The modification for the  first order moment $\overline v_i$  would be the fifth order accuracy in the smooth regions, and more detailed derivation can refer to the literature \cite{vozq}.

\textbf{Step 2.} Reconstruct the values of  the solutions $u$ at the four Gauss-Lobatto points.

We use the same stencils $S_0, S_1, S_2$ as Step 1, then,  if one of the cells in stencil $S_0$ is identified as a troubled-cell, we would reconstruct $u^\pm_{i\mp1/2}$ using the HWENO methodology in Step 2.1; otherwise we directly reconstruct $u^\pm_{i\mp1/2}$  by the linear approximation method described in Step 2.2. And the reconstruction procedure for  $u_{i\pm\sqrt5/10}$ is given in Step 2.3.

\textbf{Step 2.1.} The new HWENO reconstruction for $u^\pm_{i\mp1/2}$.

If one of the cells in stencil $S_0$ is identified as a troubled-cell, $u^\pm_{i\mp1/2}$ is reconstructed by the next  HWENO procedure. For simplicity, we only present the detailed procedure of the reconstruction for  $u^-_{i+1/2}$, while the reconstruction for $u^+_{i-1/2}$ is mirror symmetric with respect to $x_i$. Noticed that we have modified the first order moment in the troubled-cells, then, we would use these information here. We now reconstruct three polynomials $p_0(x), p_1(x), p_2(x)$ on $S_0, S_1, S_2$, respectively, satisfying
\begin{equation*}
\label{HWEp1}
\begin{split}
 \frac{1}{\Delta x} \int_{I_{i+j}} p_0(x)dx &= \overline u_{i+j}, \ \frac{1}{\Delta x} \int_{I_{i+j}} p_0(x)\frac{x-x_{i+j}} {\Delta x} dx = \overline v_{i+j},\quad  j=-1,0,1,\\
\frac{1}{\Delta x} \int_{I_{i+j}} p_1(x)dx &= \overline u_{i+j},\ j=-1,0, \quad  \frac{1}{\Delta x} \int_{I_{i}} p_1(x)\frac{x-x_{i}} {\Delta x}dx = \overline v_{i},\\
\frac{1}{\Delta x} \int_{I_{i+j}} p_2(x)dx &= \overline u_{i+j},\ j=0,1, \quad \frac{1}{\Delta x}  \int_{I_{i}} p_2(x)\frac{x-x_{i}} {\Delta x}dx = \overline v_{i}.\\
\end{split}
\end{equation*}
In terms of the above requirements, we first give the values of these polynomials at the point $x_{i+1/2}$, following as
\begin{equation*}
\begin{split}
p_0(x_{i+1/2})&=\frac {13} {108}\overline u_{i-1}+\frac 7{12}\overline u_i+\frac 8{27}\overline u_{i+1}+ \frac{25}{54}\overline v_{i-1}+\frac {241}{54}\overline v_i-\frac {28}{27}\overline v_{i+1},\\
p_1(x_{i+1/2})&=\frac16\overline u_{i-1}+\frac56\overline u_i+8 \overline v_i,\\
p_2(x_{i+1/2})&= \frac 5 6 \overline u_i+ \frac 1 {6}  \overline u_{i+1}+4 \overline v_i.\\
\end{split}
\end{equation*}
Using the next new HWENO methodology, we can use any positive linear weights satisfying $\gamma_0+ \gamma_1+\gamma_2=1$, then, we compute the smoothness indicators $\beta_n$ in the same ways, and the formula of the smoothness indicators has been given in (\ref{GHYZ}) on Step 1, then, their expressions are given as follows,
\begin{equation*}
\label{GHYZP}
\left \{
\begin{aligned}
\beta_0=& (\frac{19}{108} \overline u_{i-1}-\frac{19}{108} \overline u_{i+1} +\frac{31}{54}\overline v_{i-1}-\frac{241}{27}\overline v_{i}+\frac{31}{54}\overline v_{i+1})^2+(\frac94 \overline u_{i-1}-\frac92 \overline u_{i}+\frac94 \overline u_{i+1}+\\
&\frac{15}2\overline v_{i-1}-\frac{15}2\overline v_{i+1})^2+(\frac{70}{9} \overline u_{i-1}-\frac{70}{9} \overline u_{i+1} +\frac{200}{9}\overline v_{i-1}+\frac{1280}{9}\overline v_{i}+\frac{200}{9}\overline v_{i+1})^2+\\
&\frac{1}{12}(\frac52 \overline u_{i-1}-5\overline u_{i}+\frac52 \overline u_{i+1}+9\overline v_{i-1}-9\overline v_{i+1})^2+\frac{1}{12}(\frac{175}{18} \overline u_{i-1}-\frac{175}{18} \overline u_{i+1} +\frac{277}{9}\overline v_{i-1}+\\
&\frac{1546}{9}\overline v_{i}+\frac{277}{9}\overline v_{i+1})^2+\frac{1}{180}(\frac{95}{18} \overline u_{i-1}-\frac{95}{18}  \overline u_{i+1} +\frac{155}{9}\overline v_{i-1}+\frac{830}{9}\overline v_{i}+\frac{155}{9}\overline v_{i+1})^2+\\
&\frac{109341}{175}(\frac58 \overline u_{i-1}-\frac54\overline u_{i}+\frac58 \overline u_{i+1}+\frac{15}{4}\overline v_{i-1}-\frac{15}{4}\overline v_{i+1})^2+\frac{27553933}{1764}(\frac{35}{36} \overline u_{i-1}-\frac{35}{36}  \overline u_{i+1}+\\
&\frac{77}{18}\overline v_{i-1}+\frac{133}{9}\overline v_{i}+\frac{77}{18}\overline v_{i+1})^2,\\
\beta_1=&144\overline v_{i}^2+\frac {13}{3}(\overline u_{i-1}-\overline u_{i}+12\overline v_{i} )^2,\\
\beta_2=&144\overline v_{i}^2+\frac {13}{3}(\overline u_{i}-\overline u_{i+1}+12\overline v_{i} )^2.\\
\end{aligned}
\right.
\end{equation*}
We bring the same parameter $\tau$ to define the absolute difference between $\beta_{0}$, $\beta_{1}$ and $\beta_{2}$, and the formula is given in (\ref{tao5}), then, the nonlinear weights are computed as
\begin{equation*}
\label{90}
\omega_n=\frac{\bar\omega_n}{\sum_{\ell=0}^{2}\bar\omega_{\ell}},
\ \mbox{with} \ \bar\omega_{n}=\gamma_{n}(1+\frac{\tau}{\beta_{n}+\varepsilon}), \ n=0,1,2.
\end{equation*}
Here, $\varepsilon$ is a small positive number taken as $10^{-6}$. Finally,  the value of $u^-_{i+1/2}$ is reconstructed by
\begin{equation*}
  u^-_{i+1/2} =\omega_0 \left( \frac 1 {\gamma_0}p_0(x_{i+1/2})  - \sum_{n=1}^{2}\frac {\gamma_n} {\gamma_0} p_n(x_{i+1/2}) \right) + \sum_{n=1}^{2}\omega_n p_n(x_{i+1/2}).
\end{equation*}

\textbf{Step 2.2.} The linear approximation for $u^\mp_{i\pm1/2}$.

If neither cell in stencil $S_0$ is identified as troubled-cell,  we will use the linear approximation for $ u^\mp_{i\pm1/2}$, which means we only need to use the high degree polynomial $p_0(x)$ obtained in Step 2.1, then, we have
\begin{equation*}
u^+_{i-1/2}= p_0(x_{i-1/2})=\frac 8{27}\overline u_{i-1}+\frac 7{12}\overline u_i+\frac {13} {108}\overline u_{i+1}+\frac {28}{27} \overline v_{i-1}-\frac {241}{54}\overline v_i-\frac{25}{54}\overline v_{i+1},
\end{equation*}
and
\begin{equation*}
u^-_{i+1/2}=p_0(x_{i+1/2})=\frac {13} {108}\overline u_{i-1}+\frac 7{12}\overline u_i+\frac 8{27}\overline u_{i+1}+ \frac{25}{54}\overline v_{i-1}+\frac {241}{54}\overline v_i-\frac {28}{27}\overline v_{i+1}.
\end{equation*}

\textbf{Step 2.3.} The linear approximation for $u_{i\pm\sqrt5/10}$.

We would reconstruct $u_{i\pm\sqrt5/10}$ using the linear approximation for all cells, then, $u_{i\pm\sqrt5/10}$ are approximated by
\begin{equation*}
\begin{split}
u_{i-\sqrt5/10}= p_0(x_{i-\sqrt5/10})&=-(\frac{101}{5400}\sqrt5+\frac{1}{24})\overline u_{i-1}+\frac{13}{12}\overline u_i+(\frac{101}{5400}\sqrt5-\frac{1}{24})\overline u_{i+1}-\\
&\quad (\frac{3}{20}+\frac{841}{13500}\sqrt5)\overline v_{i-1}-\frac{10289}{6750}\sqrt5\overline v_i+(\frac{3}{20}-\frac{841}{13500}\sqrt5)\overline v_{i+1},
\end{split}
\end{equation*}
and
\begin{equation*}
\begin{split}
u_{i+\sqrt5/10}= p_0(x_{i+\sqrt5/10})&=(\frac{101}{5400}\sqrt5-\frac{1}{24})\overline u_{i-1}+\frac{13}{12}\overline u_i-(\frac{101}{5400}\sqrt5+\frac{1}{24})\overline u_{i+1}+\\
&\quad (\frac{841}{13500}\sqrt5-\frac{3}{20})\overline v_{i-1}+\frac{10289}{6750}\sqrt5\overline v_i+(\frac{3}{20}+\frac{841}{13500}\sqrt5)\overline v_{i+1}.
\end{split}
\end{equation*}

\textbf{Step 3.} Discretize the semi-discrete scheme (\ref{odeH})
in time by the third order TVD Runge-Kutta method \cite{so1}
\begin{eqnarray}
\label{RK}\left \{
\begin{array}{lll}
     u^{(1)} & = & u^n + \Delta t L(u^n),\\
     u^{(2)} & = & \frac 3 4u^n + \frac 1 4u^{(1)}+\frac 1 4\Delta t L(u^{(1)}),\\
     u^{(n+1)} & = &\frac 1 3 u^n +  \frac 2 3u^{(2)} +\frac 2 3\Delta t L(u^{(2)}).
\end{array}
\right.
\end{eqnarray}

{\bf \em Remark 1:} The KXRCF troubled-cells indicator can catch the discontinuities well. For one dimensional scalar equation, the solution $u$ is defined as the indicator variable, then $\overrightarrow v$ is $f'(u)$.  For one dimensional Euler equations, the density $\rho$ and the energy $E$ are set as the indicator variables, respectively, then $\overrightarrow v$ is the velocity $\mu$ of the fluid.

{\bf \em Remark 2:} For the systems, such as the one dimensional compressible Euler equations, all HWENO procedures are performed on the local characteristic directions to avoid the oscillations nearby discontinuities, while the linear approximation procedures are computed in each component straightforwardly.

\subsection{Two dimensional case}

We first consider two dimensional scalar hyperbolic conservation laws
  \begin{equation}
\label{EQ2} \left\{
\begin{array}
{ll}
u_t+ f(u)_x+g(u)_y=0, \\
u(x,y,0)=u_0(x,y), \\
\end{array}
\right.
\end{equation}
then, we divide the computing domain by uniform meshes $I_{i,j}$=$ [x_{i-1/2},x_{i+1/2}] \times [y_{j-1/2},y_{j+1/2}]$ for simplicity. The mesh sizes are $ \Delta x =x_{i+1/2}- x_{i-1/2}$ in the $x$ direction and  $ \Delta y =y_{j+1/2}- y_{j-1/2}$ in the $y$ direction. The cell center $(x_i,y_j)=(\frac{x_{i-1/2}+x_{i+1/2}}{2}, \frac {y_{j-1/2} +y_{j+1/2}}{2})$. $x_i+a\Delta x$ is simplified as $x_{i+a}$ and $y_j+b\Delta y$  is set  as $y_{j+b}$.

 Since the variables of the HWENO scheme are the zeroth and first order moments, we multiply the governing equation (\ref{EQ2}) by $\frac{1}{\Delta x\Delta y}$,  $\frac {x-x_i} {(\Delta x)^2\Delta y}$ and $\frac {y-y_j} {\Delta x(\Delta y)^2}$ on both sides,  respectively, then, we integrate them over $I_{i,j}$ and apply the integration by parts. In addition, we approximate the values of the flux at the points on the interface of $I_{i,j}$ by the numerical flux. Finally, the semi-discrete finite volume HWENO scheme is
\begin{equation}
\label{ode2}
\left\{
\begin{aligned}
\frac{d \overline u_{i,j}(t)}{dt}&=-\frac{1} {\Delta x \Delta y}  \int_{y_{j-1/2}}^{y_{j+1/2}}[\hat f(u(x_{i+1/2},y))-\hat f(u(x_{i-1/2},y))]dy\\
&-\frac{1} {\Delta x \Delta y}  \int_{x_{i-1/2}}^{x_{i+1/2}}[\hat g(u(x,y_{j+1/2}))-\hat g(u(x,y_{j-1/2}))]dx,\\
\frac{d \overline v_{i,j}(t)}{dt}&=-\frac{1} {2\Delta x \Delta y}  \int_{y_{j-1/2}}^{y_{j+1/2}}[\hat f(u(x_{i-1/2},y))+\hat f(u(x_{i+1/2},y))]dy+\frac 1{{\Delta x}^2 \Delta y}\int_{I_{i,j}}f(u)dxdy\\
&-\frac{1} {\Delta x \Delta y}  \int_{x_{i-1/2}}^{x_{i+1/2}}\frac{(x-x_i)}{\Delta x}[\hat g(u(x,y_{j+1/2}))-\hat g(u(x,y_{j-1/2}))]dx,\\
\frac{d \overline w_{i,j}(t)}{dt}&=-\frac{1} {\Delta x \Delta y}  \int_{y_{j-1/2}}^{y_{j+1/2}}\frac{(y-y_j)}{\Delta y}[\hat f(u(x_{i+1/2},y))-\hat f(u(x_{i-1/2},y))]dy\\
&-\frac{1} {2\Delta x \Delta y}  \int_{x_{i-1/2}}^{x_{i+1/2}}[\hat g(u(x,y_{j-1/2}))+\hat g(u(x,y_{j+1/2}))]dx+\frac 1{{\Delta x \Delta y}^2}\int_{I_{i,j}}g(u)dxdy.
\end{aligned}
\right.
\end{equation}
The initial conditions are $ \overline u_{i,j}(0) =$$ \frac 1 {\Delta x \Delta y} \int_{I_{i,j}} u_0(x,y) dxdy$, $ \overline v_{i,j}(0) =$$ \frac 1 {\Delta x \Delta y} \int_{I_{i,j}} u_0(x,y) \frac{x-x_i} { \Delta x} dxdy$ and $ \overline w_{i,j}(0) =$$ \frac 1 {\Delta x \Delta y} \int_{I_{i,j}} u_0(x,y) \frac{y-y_j} { \Delta y} dxdy$. Here, $\overline u_{i,j}(t)$ is the zeroth order moment defined as $\frac 1 {\Delta x \Delta y} $$\int_{I_{i,j}} u(x,y,t)dxdy$; $\overline v_{i,j}(t)$ and $\overline w_{i,j}(t)$ are the first order moments in the $x$ and $y$ directions taken as $\frac 1 {\Delta x \Delta y}$$\int_{I_{i,j}} u(x,y,t)$$ \frac{x-x_i} { \Delta x} dxdy$ and $\frac 1 {\Delta x \Delta y}$$\int_{I_{i,j}} u(x,y,t)$$\frac{y-y_j}{\Delta y}dxdy$, respectively.
$\hat f(u(x_{i+1/2},y))$ and $ \hat g(u(x,y_{j+1/2}))$ are the numerical flux to approximate the values of  $ f(u(x_{i+1/2},y)) $ and  $g(u(x,y_{j+1/2}))$, respectively.

Now, we approximate the integral terms of equations (\ref{ode2}) by 3-point Gaussian numerical integration. More explicitly, the integral terms are approximated by
\begin{equation*}
 \frac 1{\Delta x\Delta y} \int_{I_{i,j}}f(u)dxdy \approx \sum_{k=1}^{3}\sum_{l=1}^{3}\omega_k\omega_l f(u(x_{G_k},y_{G_l})),
\end{equation*}
\begin{equation*}
  \int_{y_{j-1/2}}^{y_{j+1/2}}\hat f(u(x_{i+1/2},y))dy\approx \Delta y \sum_{k=1}^{3} \omega_k \hat f(u(x_{i+1/2},y_{G_k})),
\end{equation*}
in which $\omega_1=\frac5{18}$, $\omega_2=\frac4{9}$ and $\omega_3=\frac5{18}$  are the quadrature weights, and the coordinates of the Gaussian points are
\begin{equation*}
x_{G_1}=x_{i-\frac{\sqrt{15}}{10}},\ x_{G_2}=x_{i},\ x_{G_3}=x_{i+\frac{\sqrt{15}}{10}}; \quad y_{G_1}=y_{j-\frac{\sqrt{15}}{10}},\ y_{G_2}=y_{j},\ y_{G_3}=y_{j+\frac{\sqrt{15}}{10}}.
\end{equation*}

The numerical fluxes at the interface points in each directions are approximated by the Lax-Friedrichs method:
\begin{equation*}\label{2dflx}
  \hat f(u(G_b))=\frac 1 2[f(u^-(G_b))+f(u^+(G_b))]-\frac{\alpha}2(u^+(G_b)-u^-(G_b)),
\end{equation*}
and
\begin{equation*}\label{2dfly}
\hat g(u(G_b))=\frac 1 2[g(u^-(G_b))+g(u^+(G_b))]-\frac{\beta}2(u^+(G_b)-u^-(G_b)).
\end{equation*}
Here, $\alpha= \max_u|f'(u)|$, $\beta= \max_u|g'(u)|$, and $G_b$ is the Gaussian point on the interface of the cell $I_{i,j}$.

Now, we first present the detailed spatial reconstruction  for the semi-discrete scheme  (\ref{ode2}) in Steps 4 and 5, then, we introduce the methodology of time discretization in Step 6.

\textbf{Step 4.} Identify the troubled-cell  and  modify the first order moments in the troubled-cell.

 We also use the KXRCF troubled-cell indicator \cite{LJJN} to  identify the discontinuities, and the detailed implementation procedures for two dimensional problems had been introduced in the hybrid HWENO scheme \cite{ZCQHH}.

 If the cell $I_{i,j}$ is identified as a troubled-cell, we would modify the first order moments $\overline v_{i,j}$ and $\overline w_{i,j}$. We can  modify the first order moments employing dimensional by dimensional manner. For example, we use these information $\overline u_{i-1,j}$, $\overline u_{i,j}$, $\overline u_{i+1,j}$, $\overline v_{i-1,j}$, $\overline v_{i+1,j}$ to  modify $\overline v_{i,j}$, but employ $\overline u_{i,j-1}$, $\overline u_{i,j}$, $\overline u_{i,j+1}$, $\overline w_{i,j-1}$, $\overline w_{i,j+1}$ to  reconstruct $\overline w_{i,j}$, and the procedures are the same as one dimensional case.

\textbf{Step 5.} Reconstruct the point values of the solutions $u$ at the Gaussian points.

Based on the formula of the semi-discrete scheme  (\ref{ode2}), it means that we need to reconstruct the point values of $u^\pm(x_{i\mp1/2},y_{G_{1,2,3}})$, $u^\pm(x_{G_{1,2,3}},y_{j\mp1/2})$ and $u(x_{G_{1,2,3}},y_{G_{1,2,3}})$ in the cell $I_{i,j}$. If one of the cells in the big stencil is identified as a troubled-cell in Step 4, we would reconstruct the points values of solutions $u$ at the interface points of the cell $I_{i,j}$ by the HWENO methodology in Step 5.1; otherwise we directly use linear approximation at these interface points in Step 5.2. And we employ linear approximation straightforwardly for internal reconstructed points introduced in Step 5.3.

\textbf{Step 5.1.} Reconstruct the point values of the solutions $u$ at the interface points by a new HWENO methodology.

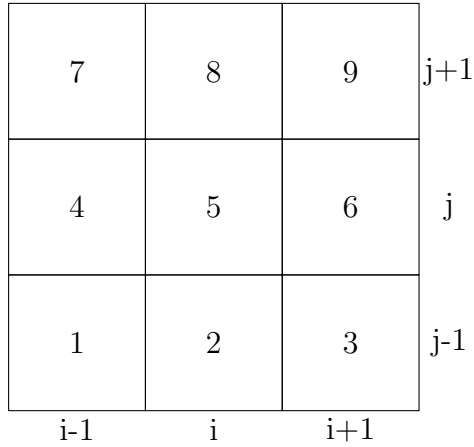
\begin{figure}
\tikzset{global scale/.style={
 scale=#1,every node/.append style={scale=#1}}}
  \centering
   \begin{tikzpicture}[global scale = 1]
     \draw(0,0)rectangle+(1.8,1.8);\draw(1.8,0)rectangle+(1.8,1.8);\draw(3.6,0)rectangle+(1.8,1.8);
     \draw(0,1.8)rectangle+(1.8,1.8);\draw(1.8,1.8)rectangle+(1.8,1.8);\draw(3.6,1.8)rectangle+(1.8,1.8);
     \draw(0,3.6)rectangle+(1.8,1.8);\draw(1.8,3.6)rectangle+(1.8,1.8);\draw(3.6,3.6)rectangle+(1.8,1.8);
     \draw(0.9,0.9)node{1};\draw(2.7,0.9)node{2};\draw(4.5,0.9)node{3};
     \draw(0.9,2.7)node{4};\draw(2.7,2.7)node{5};\draw(4.5,2.7)node{6};
     \draw(0.9,4.5)node{7};\draw(2.7,4.5)node{8};\draw(4.5,4.5)node{9};
     \draw(0.9,-0.25)node{i-1};\draw(2.7,-0.25)node{i};\draw(4.5,-0.25)node{i+1};
     \draw(5.8,0.9)node{j-1};\draw(5.8,2.7)node{j};\draw(5.8,4.5)node{j+1};
   \end{tikzpicture}
 \caption{The big stencil $S_0$ and its new labels.}
 \label{2dbig}
\end{figure}

\begin{figure}
\tikzset{global scale/.style={
 scale=#1,every node/.append style={scale=#1}}}
  \centering
   \begin{tikzpicture}[global scale = 1]
   \draw(0+0,0+9.6)rectangle+(1.8,1.8);\draw(1.8+0,0+9.6)rectangle+(1.8,1.8);
   \draw(0+0,1.8+9.6)rectangle+(1.8,1.8);\draw(1.8+0,1.8+9.6)rectangle+(1.8,1.8);
   \draw(0.9+0,0.9+9.6)node{1};\draw(2.7+0,0.9+9.6)node{2};
   \draw(0.9+0,2.7+9.6)node{4};\draw(2.7+0,2.7+9.6)node{5};
    \draw(0.9+0,-0.25+9.6)node{i-1};\draw(2.7+0,-0.25+9.6)node{i};
    \draw(4.0+0,0.9+9.6)node{j-1};\draw(4.0+0,2.7+9.6)node{j};
   \draw(0+5.0,0+9.6)rectangle+(1.8,1.8);\draw(1.8+5.0,0+9.6)rectangle+(1.8,1.8);
   \draw(0+5.0,1.8+0+9.6)rectangle+(1.8,1.8);\draw(1.8+5.0,1.8+9.6)rectangle+(1.8,1.8);
   \draw(0.9+5.0,0.9+9.6)node{2};\draw(2.7+5.0,0.9+9.6)node{3};
   \draw(0.9+5.0,2.7+9.6)node{5};\draw(2.7+5.0,2.7+9.6)node{6};
   \draw(0.9+5.0,-0.25+9.6)node{i};\draw(2.7+5.0,-0.25+9.6)node{i+1};
   \draw(4.0+5.0,0.9+9.6)node{j-1};\draw(4.0+5.0,2.7+9.6)node{j};

     \draw(0+0,0+14.4)rectangle+(1.8,1.8);\draw(1.8+0,0+14.4)rectangle+(1.8,1.8);
     \draw(1.8+0,1.8+14.4)rectangle+(1.8,1.8);\draw(0+0,1.8+14.4)rectangle+(1.8,1.8);
    \draw(0.9+0,0.9+14.4)node{4};\draw(2.7+0,0.9+14.4)node{5};
  \draw(0.9+0,2.7+14.4)node{7};\draw(2.7+0,2.7+14.4)node{8};
     \draw(0.9+0,-0.25+14.4)node{i-1};\draw(2.7+0,-0.25+14.4)node{i};
    \draw(4+0,0.9+14.4)node{j};\draw(4+0,2.7+14.4)node{j+1};
   \draw(0+5.0,0+14.4)rectangle+(1.8,1.8);\draw(1.8+5.0,0+14.4)rectangle+(1.8,1.8);
   \draw(0+5.0,1.8+14.4)rectangle+(1.8,1.8);\draw(1.8+5.0,1.8+14.4)rectangle+(1.8,1.8);
    \draw(0.9+5.0,0.9+14.4)node{5}; \draw(2.7+5.0,0.9+14.4)node{6};
    \draw(0.9+5.0,2.7+14.4)node{8}; \draw(2.7+5.0,2.7+14.4)node{9};
     \draw(0.9+5.0,-0.25+14.4)node{i};\draw(2.7+5.0,-0.25+14.4)node{i+1};
    \draw(4+5.0,0.9+14.4)node{j};\draw(4+5.0,2.7+14.4)node{j+1};
   \end{tikzpicture}
 \caption{The four small stencils and these respective labels. From left to right and bottom to top are the stencils: $S_1,...,S_4$.}
 \label{2dsmall}
\end{figure}
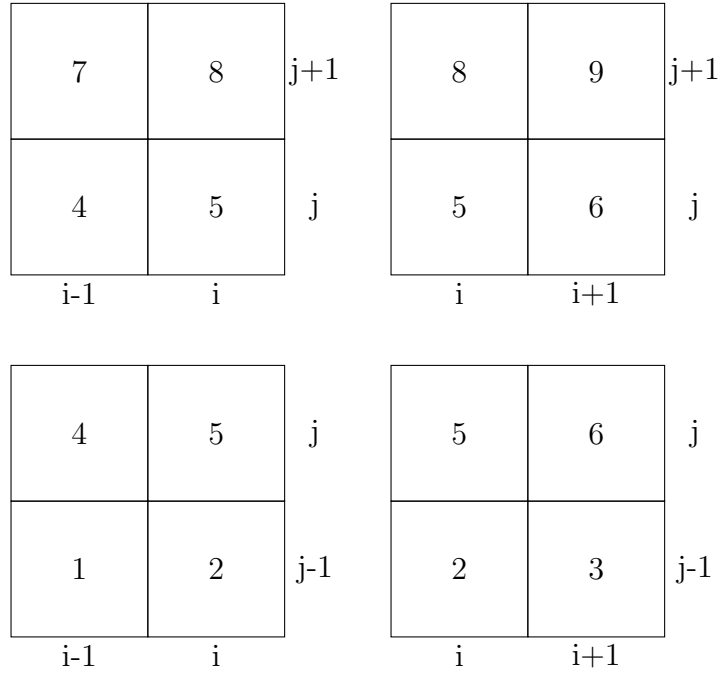

If one of the cells in big stencil is identified as a troubled-cell,  the points values of solutions $u$ at the interface points of the cell $I_{i,j}$ are reconstructed by the next new HWENO methodology.
We first give the big stencil $S_0$ in Figure \ref{2dbig}, and we rebel the cell $I_{i,j}$ and its neighboring cells  as $I_1,...,I_9$ for simplicity. Particularly, the new label of the cell $I_{i,j}$ is $I_5$. In the next procedures, we take $G_k$ to represent the specific points where we want to reconstruct.  We also give four small stencils $S_1,...,S_4$ shown in Figure \ref{2dsmall}. Noticed that we only use five candidate stencils, but the hybrid HWENO scheme \cite{ZCQHH} needed to use eight small stencils. Now, we  construct a  quartic reconstruction polynomial $p_0(x,y)$ $\in span$ $\{1,x,y,x^2,xy,y^2,x^3,x^2y,xy^2,y^3,x^4,x^3y,x^2y^2,xy^3, y^4\}$ on the big stencil $S_0$ and four quadratic polynomials $p_1(x,y),...,p_4(x,y)$ $\in span\{1,x,y,x^2,xy,y^2\}$ on the four small stencils $S_1,...,S_4$, respectively. These polynomials satisfy the following conditions:
\begin{equation*}
\begin{array}{ll}
\frac{1}{\Delta x \Delta y}\int_{I_k}p_n(x,y)dxdy=\overline u_k,            \\
 \frac{1}{\Delta x \Delta y}\int_{I_{k_x}}p_n(x,y)\frac{(x-x_{k_x})}{\Delta x}dxdy=\overline v_{k_x}, \quad \frac{1}{\Delta x \Delta y}\int_{I_{k_y}}p_n(x,y)\frac{(y-y_{k_y})}{\Delta y}dxdy=\overline w_{k_y}, \\
\end{array}
\end{equation*}
for
  \begin{equation*}
  \begin{array}{ll}
n=0,\quad k=1,...,9,\ k_x=k_y=2,4,5,6,8;\\
n=1,\quad k=1,2,4,5,\ k_x=k_y=5; \quad n=2,\quad k=2,3,5,6,\ k_x=k_y=5;\\
n=3,\quad k=4,5,7,8,\ k_x=k_y=5; \quad n=4,\quad k=5,6,8,9,\ k_x=k_y=5.
\end{array}
\end{equation*}
For the quartic polynomial $p_0(x,y)$, we can obtain it by requiring that it matches the zeroth order moments on the cell $I_1,...,I_9$, the first order moments on the cell $I_5$ and others are in a least square sense \cite{hs}. For the four quadratic polynomials, we can directly obtain the expressions of $p_n(x,y)$ $(n=1,...,4)$ by the above corresponding requirements, respectively.

Similarly as in the one dimensional case, the new HWENO method can use any artificial positive linear weights (the sum equals 1), while the hybrid HWENO scheme \cite{ZCQHH} needed to calculate the linear weights for 12 points using 8 small stencils determined  by a least square methodology, and the linear weights were not easy to be obtained especially for high dimensional problems or unstructured meshes. In addition, it only had the fourth order accuracy in two dimension, but the new HWENO methodology can achieve the fifth order numerical accuracy. Next, to measure how smooth the function $p_n(x, y)$ in the target cell $I_{i,j}$, we compute the smoothness indicators $\beta_n$ as the same way listed by \cite{hs}, following as
\begin{equation}\label{2dGHYZ}
  \beta_n= \sum_{|l|=1}^r|I_{i,j}|^{|l|-1} \int_{I_{i,j}}\left( \frac {\partial^{|l|}}{\partial x^{l_1}\partial y^{l_2}}p_n(x,y)\right)^2 dxdy, \quad n=0,...,4,
\end{equation}
where $l=(l_1,l_2)$, $|l|=l_1+l_2$ and $r$ is the degree of $p_n(x, y)$. Similarly, we  bring a new parameter $\tau$ to define the overall difference between $\beta_{l}$, $l=0,...,4$ as
\begin{equation}
\label{tao4}
\tau=\left(\frac{|\beta_{0}-\beta_{1}|+|\beta_{0}-\beta_{2}|+|\beta_{0}-\beta_{3}|+|\beta_{0}-\beta_{4}|}{4}\right)^2,
\end{equation}
then, the nonlinear weights are defined as
\begin{equation}
\label{99}
\omega_n=\frac{\bar\omega_n}{\sum_{\ell=0}^{4}\bar\omega_{\ell}},
\ \mbox{with} \ \bar\omega_{n}=\gamma_{n}(1+\frac{\tau}{\beta_{n}+\varepsilon}), \ n=0,...,4,
\end{equation}
in which $\varepsilon$ is taken as $10^{-6}$.  The final reconstruction of
the solutions $u$ at the interface  point $G_k$ is
\begin{equation*}
  u^*(G_k)  =\omega_0 \left( \frac 1 {\gamma_0}p_0(G_k) - \sum_{n=1}^{4}\frac {\gamma_n} {\gamma_0} p_n(G_k) \right) + \sum_{n=1}^{4}\omega_n p_n(G_k).
\end{equation*}
where "*" is "+" when $G_k$ is located on the left or bottom interface of the cell $I_{i,j}$, while "*" is "-" on the right or top interface of $I_{i,j}$.

\textbf{Step 5.2.} Reconstruct the point values of the solutions $u$ at the interface points using linear approximation.

If neither cell in the big stencil $S_0$ is  identified as a troubled-cell, the point value of the solution $u$ at the interface point $G_k$ is directly approximated by $p_0(G_k)$, and we use the same polynomial $p_0(x,y)$ given in Step 5.1.

\textbf{Step 5.3.} Reconstruct the point values of the solutions $u$ at the internal points by linear approximation straightforwardly.

We would use linear approximation for the point values of the solutions $u$ at the internal points in all cells, then, we directly employ the same  quartic polynomial $p_0(x,y)$ obtained in  Step 5.1 to approximate these point values.

\textbf{Step 6.} Discretize the  semi-discrete scheme (\ref{ode2}) in time by the third order TVD Runge-Kutta method \cite{so1}.

The  semi-discrete scheme (\ref{ode2}) is discretized by the third order TVD Runge-Kutta method in time, and the formula is given in (\ref{RK}) for the one dimensional case.

{\bf \em Remark 3:} The KXRCF indicator is suitable  for two dimensional hyperbolic conservation laws. For two dimensional scalar equation, the solution $u$ is the indicator variable. $\overrightarrow{v}$ is set as $f'(u)$ in the $x$ direction, while it is taken as $g'(u)$ in the $y$ direction. For two dimensional Euler equations, the density $\rho$ and the energy $E$ are defined as the indicator variables, respectively. $\overrightarrow{v}$ is the velocity $\mu$ in the $x$ direction of the fluid, while it is the  velocity $\nu$  in the $y$ direction.

{\bf \em Remark 4:} For the systems, such as the two dimensional compressible Euler equations, all HWENO reconstruction procedures are performed on  the  local characteristic decompositions, while linear approximation procedures are performed on component by component.

\section{Numerical tests}
\label{sec3}
\setcounter{equation}{0}
\setcounter{figure}{0}
\setcounter{table}{0}

In this section, we present the numerical results of the new hybrid HWENO scheme which is described in Section 2.  In order to fully assess the influence of the modification of the first order moment  upon accuracy,  all cells are marked as troubled-cells in Step 1 and Step 4 for one and two dimensional cases, respectively, and we denote this method as  New HWENO scheme.  We also denote HWENO scheme and the hybrid HWENO scheme which are presented in \cite{ZCQHH}.  The CFL number is set as 0.6 expect for the hybrid HWENO scheme in  the two dimensional  non-smooth tests.

\subsection{Accuracy tests}

We will present the results of HWENO, New HWENO, Hybrid HWENO and New hybrid HWENO schemes in the one and two dimensional accuracy tests. In addition, to evaluate whether the choice of the linear weights would affect the order of the new HWENO methodology or not, we use random positive linear weights (the sum equals one) at each time step for New HWENO and New hybrid HWENO schemes.

\begin{figure}
 \centerline{
\psfig{file=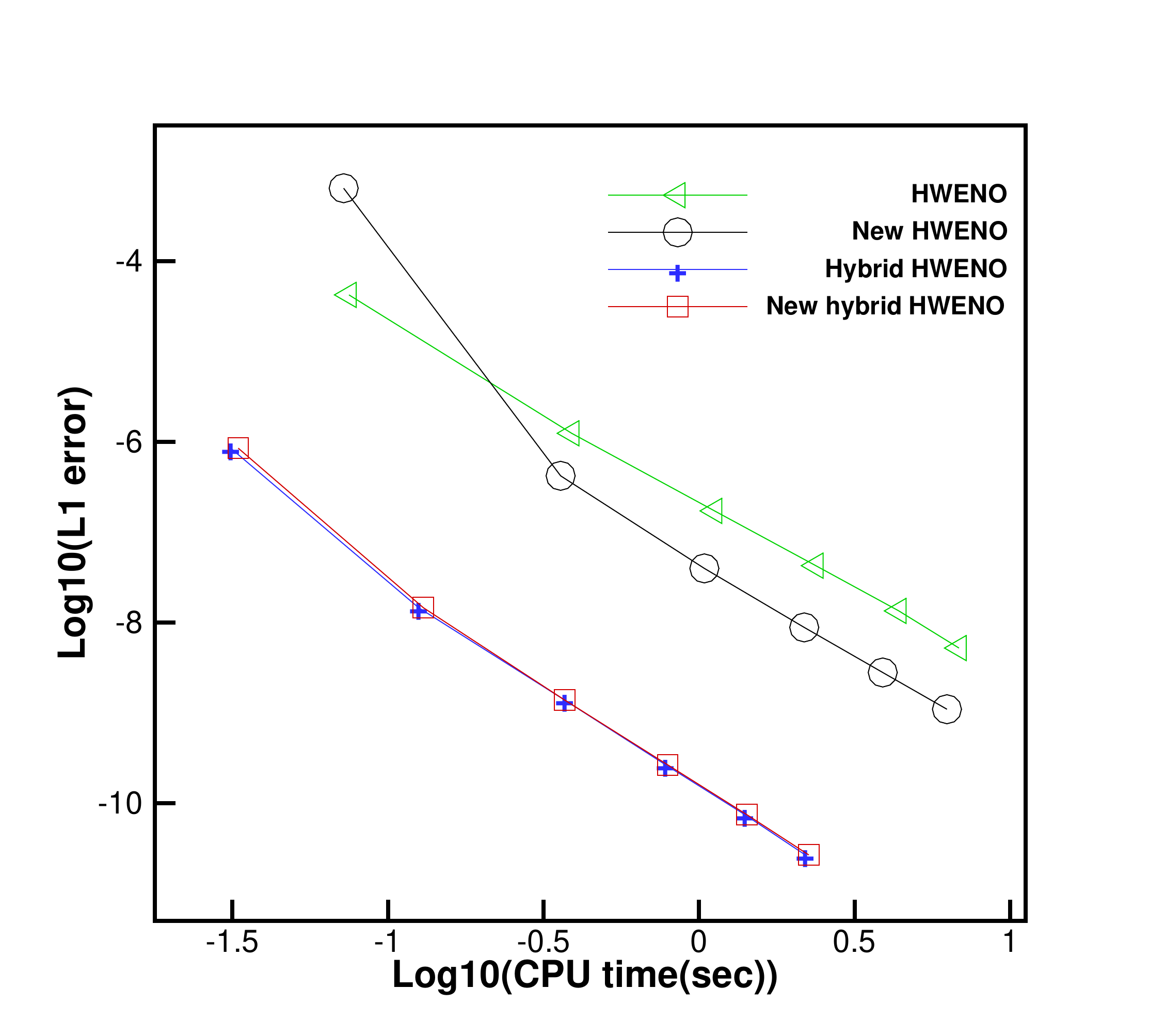,width=2.5 in}
\psfig{file=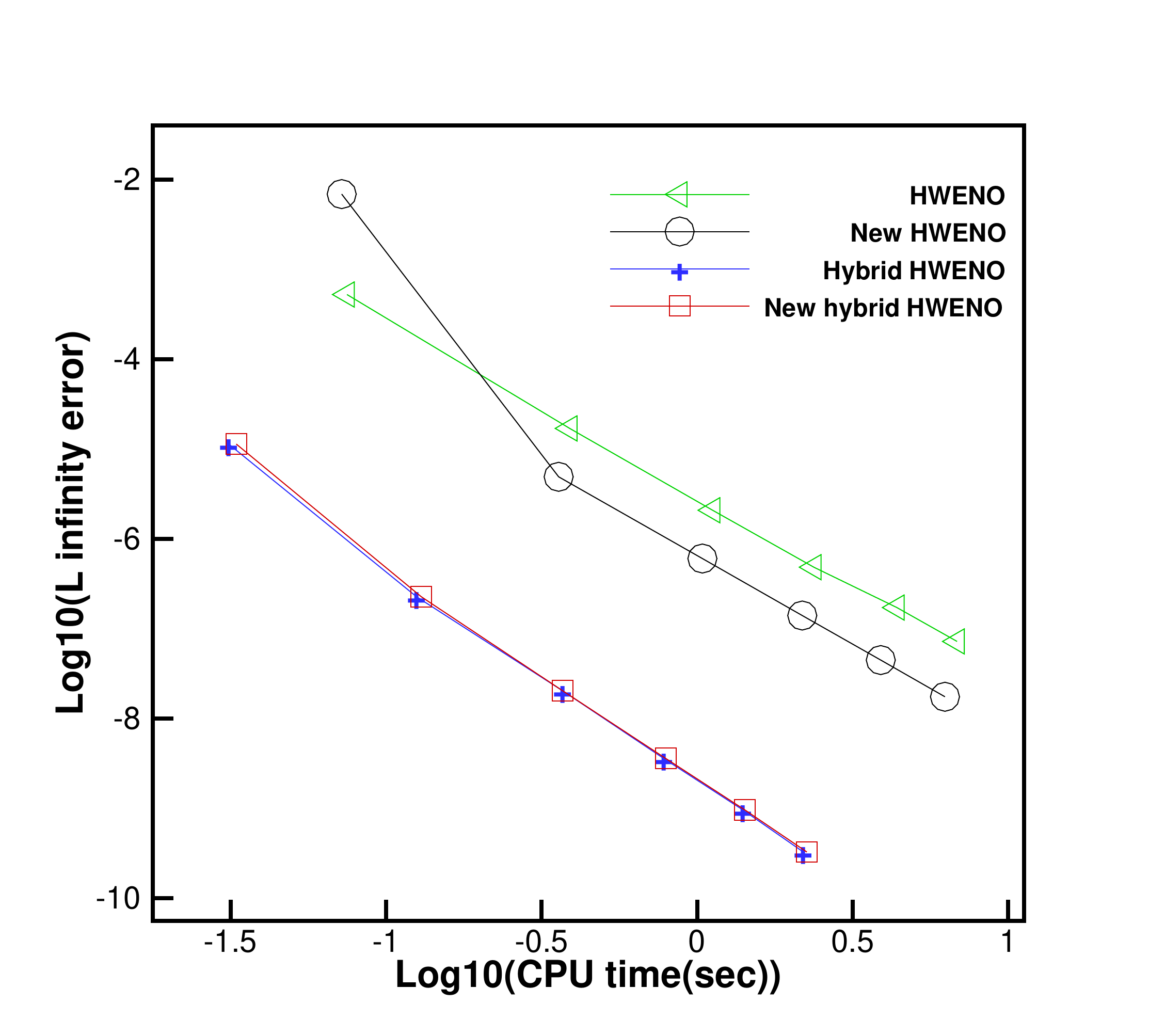,width=2.5 in}}
 \caption{1D-Burgers' equation: initial data
$u(x,0)=0.5+sin(\pi x)$. $T=0.5/\pi$. Computing times and errors. Triangle signs and a green solid  line: the results of HWENO scheme; circle signs and a black solid  line: the results of New HWENO scheme; plus signs and a blue solid  line: the results of Hybrid HWENO scheme; rectangle signs and a red solid  line: the results of New hybrid HWENO scheme.}
\label{Fburges1d_smooth}
\end{figure}

\noindent{\bf Example 3.1.} We solve the following scalar Burgers' equation:
\begin{equation}\label{1dbugers}
  u_t+(\frac {u^2} 2)_x=0, \quad 0<x<2.
\end{equation}
The initial  condition is $u(x,0)=0.5+\sin(\pi x)$ with periodic boundary condition. The computing time is $t=0.5/\pi$, in which the solution is still smooth. We give the numerical errors and orders in Table \ref{tburgers1d} with $N$ uniform meshes for HWENO, New HWENO, Hybrid HWENO and New hybrid HWENO schemes.
At first, we know that Hybrid HWENO and New hybrid HWENO schemes have same results for there are not cells which are identified as troubled-cells, therefore, they both directly use linear approximation for the spatial reconstruction. Although these HWENO schemes all have the designed fifth order accuracy, the hybrid schemes have better numerical performance with less numerical errors than the corresponding HWENO schemes, meanwhile, we can see that New HWENO scheme has less numerical errors than HWENO scheme starting with 80 meshes, which illustrates the new HWENO methodology has better numerical performance than the original HWENO method. In addition, the choice of the linear weights would not affect the order of the new HWENO methodology. Finally, we show numerical errors against CPU times by these HWENO schemes in Figure \ref{Fburges1d_smooth}, which shows two hybrid HWENO schemes have much higher efficiency than other HWENO schemes, and New HWENO scheme also has higher efficiency than HWENO scheme.
\begin{table}
\begin{center}
\caption{1D-Burgers' equation: initial data
$u(x,0)=0.5+sin(\pi x)$. HWENO schemes. $T=0.5/\pi$. $L^1$ and $L^\infty$ errors and orders. Uniform meshes with $N$ cells. }
\medskip
\begin{tabular} {lllllllll} \hline
 $N$ cells & \multicolumn{4}{l}{HWENO  scheme} & \multicolumn{4}{l}{New HWENO scheme}\\
\cline{2-5} \cline{6-9}
  &$ L^1$ error &  order & $L^\infty$error &  order &$ L^1$ error &  order & $L^\infty$ error &order\\   \hline
          40 &     4.23E-05 &     &     5.25E-04 &       							&     6.42E-04 &     &     6.89E-03 &      \\
     80 &     1.24E-06 &       5.09 &     1.70E-05 &       4.95								   &     4.20E-07 &      10.58 &     4.91E-06 &      10.45 \\
    120 &     1.72E-07 &       4.88 &     2.08E-06 &       5.17								   &     3.97E-08 &       5.82 &     6.04E-07 &       5.17 \\
    160 &     4.26E-08 &       4.85 &     4.84E-07 &       5.08								   &     8.83E-09 &       5.23 &     1.40E-07 &       5.08 \\
    200 &     1.34E-08 &       5.17 &     1.72E-07 &       4.64								   &     2.80E-09 &       5.15 &     4.47E-08 &       5.12 \\
    240 &     5.21E-09 &       5.20 &     7.22E-08 &       4.76								   &     1.10E-09 &       5.14 &     1.75E-08 &       5.16 \\
\hline
  $N$ cells & \multicolumn{4}{l}{Hybrid HWENO scheme} & \multicolumn{4}{l}{New Hybrid HWENO scheme}\\
\cline{2-5} \cline{6-9}
&$ L^1$ error &  order & $L^\infty$error &  order &$ L^1$ error &  order & $L^\infty$ error &order\\   \hline
         40 &     8.51E-07 &       &     1.14E-05 &       	&     8.51E-07 &      &     1.14E-05 &     \\
     80 &     1.46E-08 &       5.87 &     2.26E-07 &       5.65								&     1.46E-08 &       5.87 &     2.26E-07 &       5.65 \\
    120 &     1.39E-09 &       5.80 &     2.04E-08 &       5.94								&     1.39E-09 &       5.80 &     2.04E-08 &       5.94 \\
    160 &     2.66E-10 &       5.75 &     3.59E-09 &       6.03								&     2.66E-10 &       5.75 &     3.59E-09 &       6.03 \\
    200 &     7.46E-11 &       5.70 &     9.58E-10 &       5.92								&     7.46E-11 &       5.70 &     9.58E-10 &       5.92 \\
    240 &     2.68E-11 &       5.62 &     3.27E-10 &       5.90								&     2.68E-11 &       5.62 &     3.27E-10 &       5.90 \\
\hline
\end{tabular}
\label{tburgers1d}
\end{center}
\end{table}
\smallskip

\begin{figure}
 \centerline{
\psfig{file=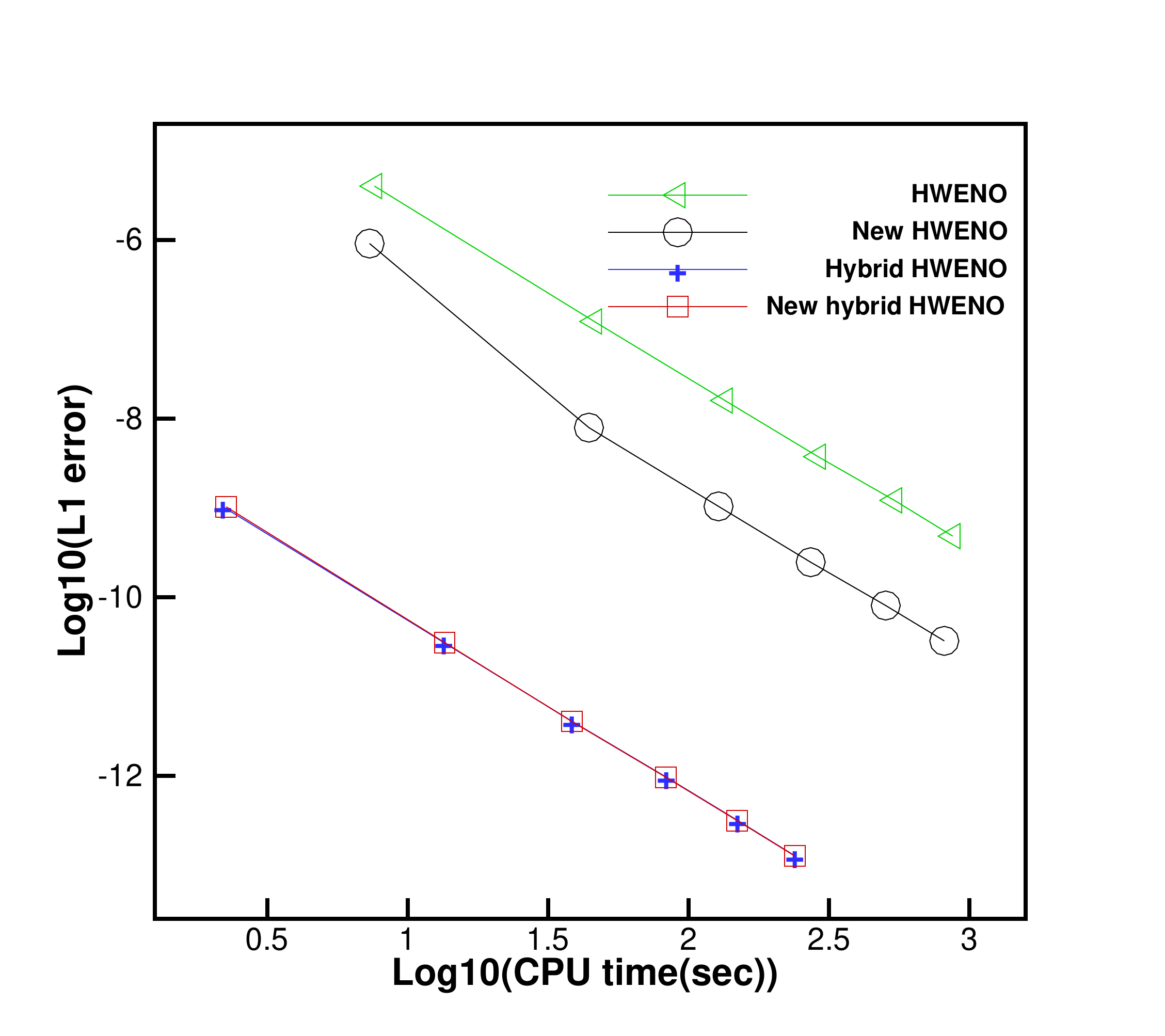,width=2.5 in}
\psfig{file=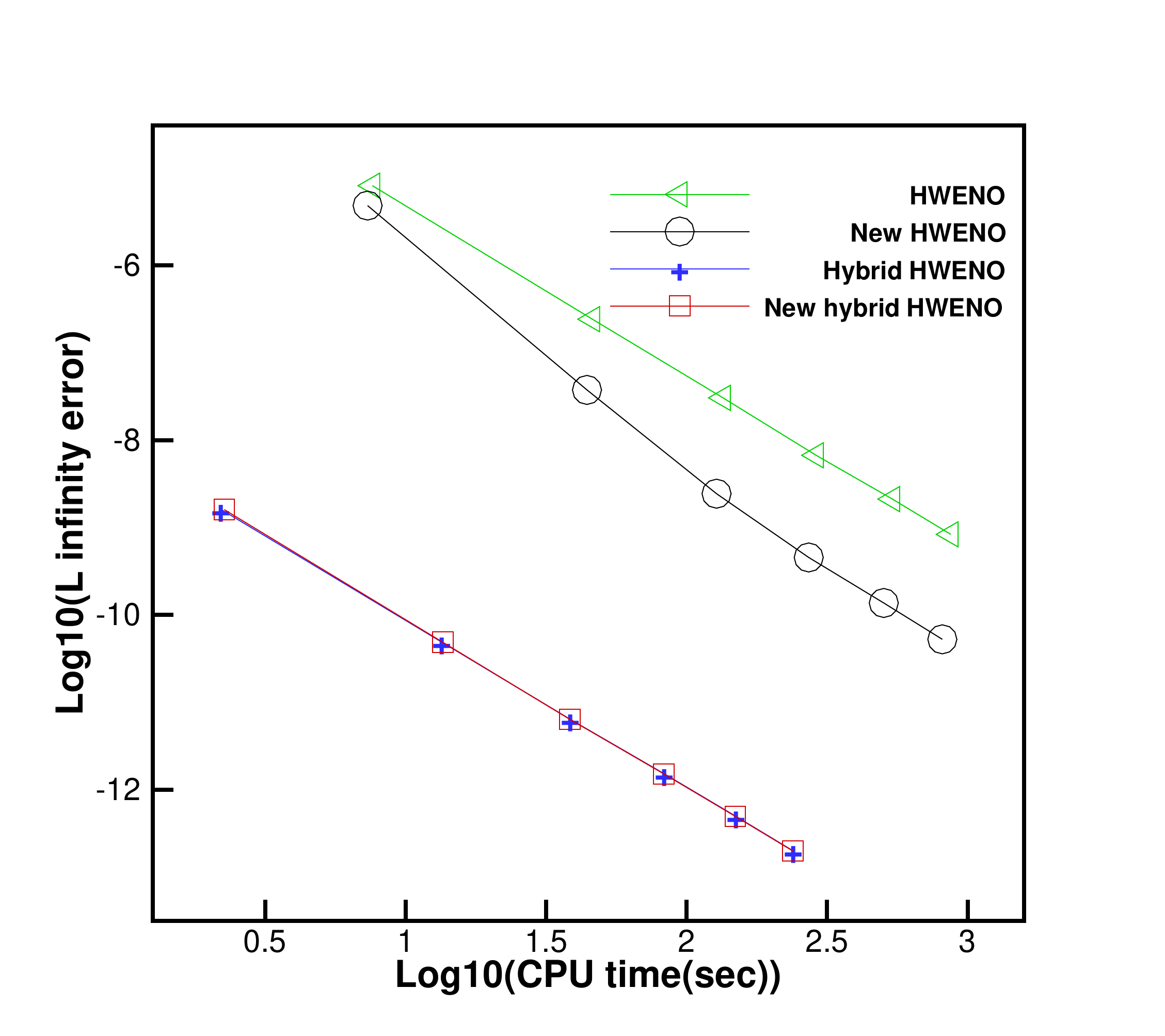,width=2.5 in}}
 \caption{1D-Euler equations: initial data
$\rho(x,0)=1+0.2\sin(\pi x)$, $\mu(x,0)=1$ and $p(x,0)=1$. $T=2$. Computing times and errors. Triangle signs and a green solid  line: the results of HWENO scheme; circle signs and a black solid  line: the results of New HWENO scheme; plus signs and a blue solid  line: the results of Hybrid HWENO scheme; rectangle signs and a red solid  line: the results of New hybrid HWENO scheme.}
\label{FEuler1d_smooth}
\end{figure}

\noindent{\bf Example 3.2.} One dimensional Euler equations:
\begin{equation}
\label{euler1}
 \frac{\partial}{\partial t}
 \left(
 \begin{array}{c}
 \rho \\
 \rho \mu \\
 E
 \end{array} \right )
 +
 \frac{\partial}{\partial x}
 \left (
\begin{array}{c}
\rho \mu \\
\rho \mu^{2}+p \\
\mu(E+p)
\end{array}
 \right )
=0,
\end{equation}
 where $\rho$ is density, $\mu$ is velocity, $E$ is total energy and $p$ is pressure. The initial conditions are $\rho(x,0)=1+0.2\sin(\pi x)$, $\mu(x,0)=1$, $p(x,0)=1$ and  $\gamma=1.4$ with periodic boundary condition. The computing domain is $ x \in [0, 2\pi]$. The exact solution is $\rho(x,t)=1+0.2\sin(\pi(x-t))$, $\mu(x,0)=1$, $p(x,0)=1$, and the computing time is up to $T=2$. We present the numerical errors and  orders of the density for the HWENO schemes in Table \ref{tEluer1d}, then, we first can see these HWENO schemes achieve the fifth order accuracy, and two hybrid HWENO schemes have same performance as they both directly use linear approximation for the spatial reconstruction, meanwhile, the hybrid schemes have less numerical errors than the corresponding HWENO schemes. In addition, New HWENO scheme has less  errors than HWENO scheme, which shows the new HWENO methodology has better performance than the original HWENO method, and random positive linear weights at each time step would not affect the order accuracy of New HWENO scheme. Finally, we give the numerical errors against CPU times by these HWENO schemes in Figure \ref{FEuler1d_smooth}, which shows Hybrid HWENO schemes have much higher efficiency with smaller numerical errors and less CPU times than other HWENO schemes, and we can see  New HWENO scheme has higher efficiency with smaller errors than HWENO scheme.
\begin{table}
\begin{center}
\caption{1D-Euler equations: initial data
$\rho(x,0)=1+0.2\sin(\pi x)$, $\mu(x,0)=1$ and $p(x,0)=1$. HWENO schemes. $T=2$. $L^1$ and $L^\infty$ errors and orders. Uniform meshes with $N$ cells. }
\medskip
\begin{tabular} {lllllllll} \hline
 $N$ cells & \multicolumn{4}{l}{HWENO  scheme} & \multicolumn{4}{l}{New HWENO scheme}\\
\cline{2-5} \cline{6-9}
&$ L^1$ error &  order & $L^\infty$error &  order &$ L^1$ error &  order & $L^\infty$ error &order\\   \hline
       40 &     4.00E-06 &  &     8.18E-06 &       						      &     9.09E-07 &  &     4.85E-06 &       \\
     80 &     1.22E-07 &       5.04 &     2.43E-07 &       5.08 	&     7.89E-09 &       6.85 &     3.76E-08 &       7.01 \\
    120 &     1.59E-08 &       5.03 &     3.05E-08 &       5.11 								      &     1.04E-09 &       5.01 &     2.44E-09 &       6.75 \\
    160 &     3.73E-09 &       5.03 &     6.71E-09 &       5.26 		      &     2.46E-10 &       5.00 &     4.54E-10 &       5.84 \\
    200 &     1.21E-09 &       5.04 &     2.12E-09 &       5.17 								      &     8.05E-11 &       5.00 &     1.37E-10 &       5.37 \\
    240 &     4.82E-10 &       5.06 &     8.35E-10 &       5.10 		      &     3.23E-11 &       5.00 &     5.25E-11 &       5.26 \\\hline
 $N$ cells & \multicolumn{4}{l}{Hybrid HWENO  scheme} & \multicolumn{4}{l}{New hybrid HWENO scheme}\\
\cline{2-5} \cline{6-9}
   &$ L^1$ error &  order & $L^\infty$error &  order &$ L^1$ error &  order & $L^\infty$ error &order\\ \hline
       40 &     1.02E-09 &  &     1.60E-09 &       						  &     1.02E-09 &      &     1.60E-09 &       \\
     80 &     3.10E-11 &       5.05 &     4.86E-11 &       5.04		  &     3.10E-11 &       5.05 &     4.86E-11 &       5.04 \\
    120 &     4.06E-12 &       5.01 &     6.37E-12 &       5.01	  &     4.06E-12 &       5.01 &     6.37E-12 &       5.01 \\
    160 &     9.61E-13 &  5.01 &     1.51E-12 &       5.01& 9.61E-13 & 5.01 &  1.51E-12 &       5.01 \\
    200 &   3.15E-13 &   5.00 &  4.94E-13 &    5.00&   3.15E-13 &   5.00 &  4.94E-13 &       5.00 \\
    240 &     1.26E-13 &       5.00 &     1.98E-13 &       5.00								  &     1.26E-13 &   5.00 &  1.98E-13 &   5.00 \\ \hline
\end{tabular}
\label{tEluer1d}
\end{center}
\end{table}
\smallskip

\begin{figure}
 \centerline{
\psfig{file=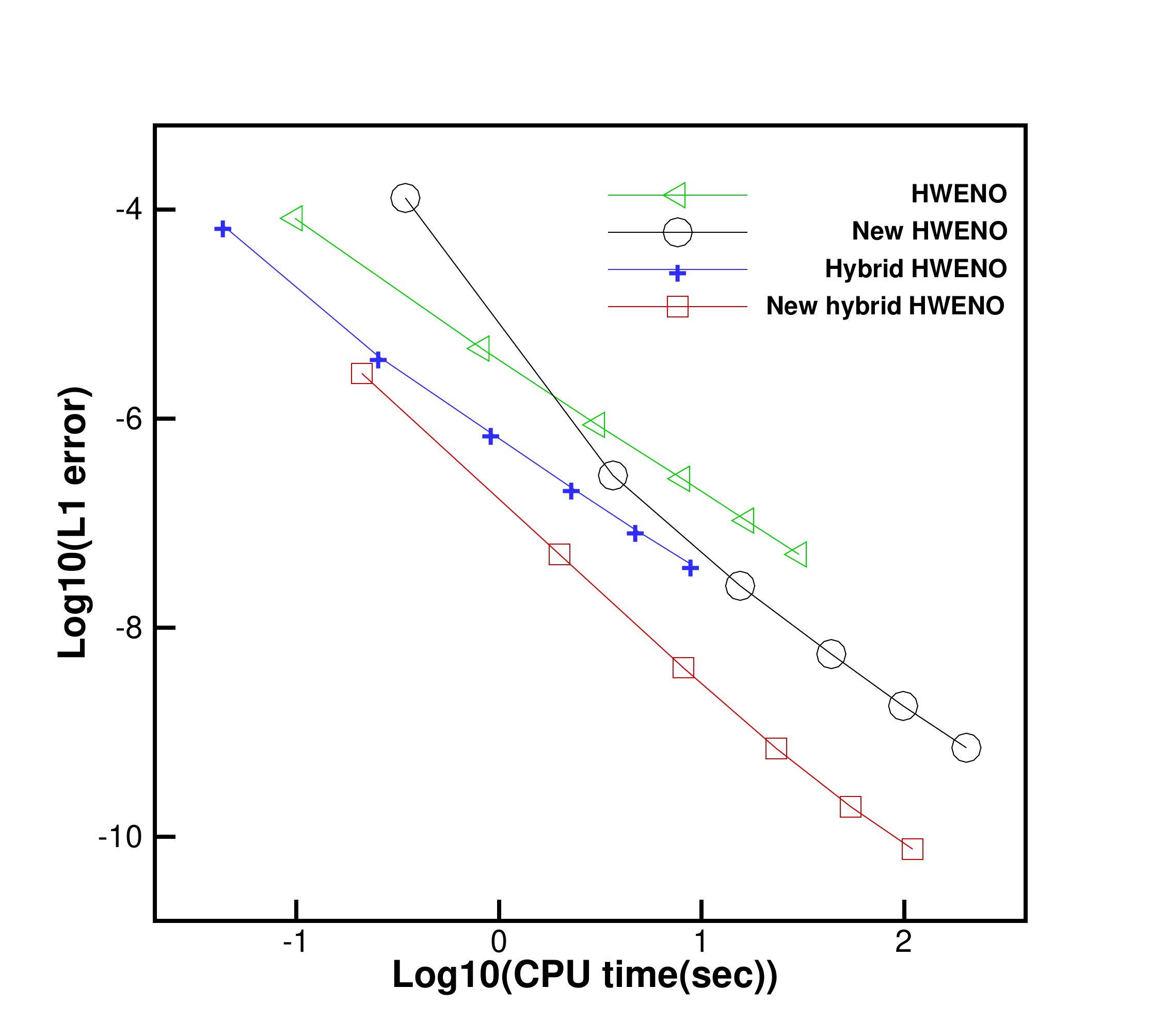,width=2.5 in}
\psfig{file=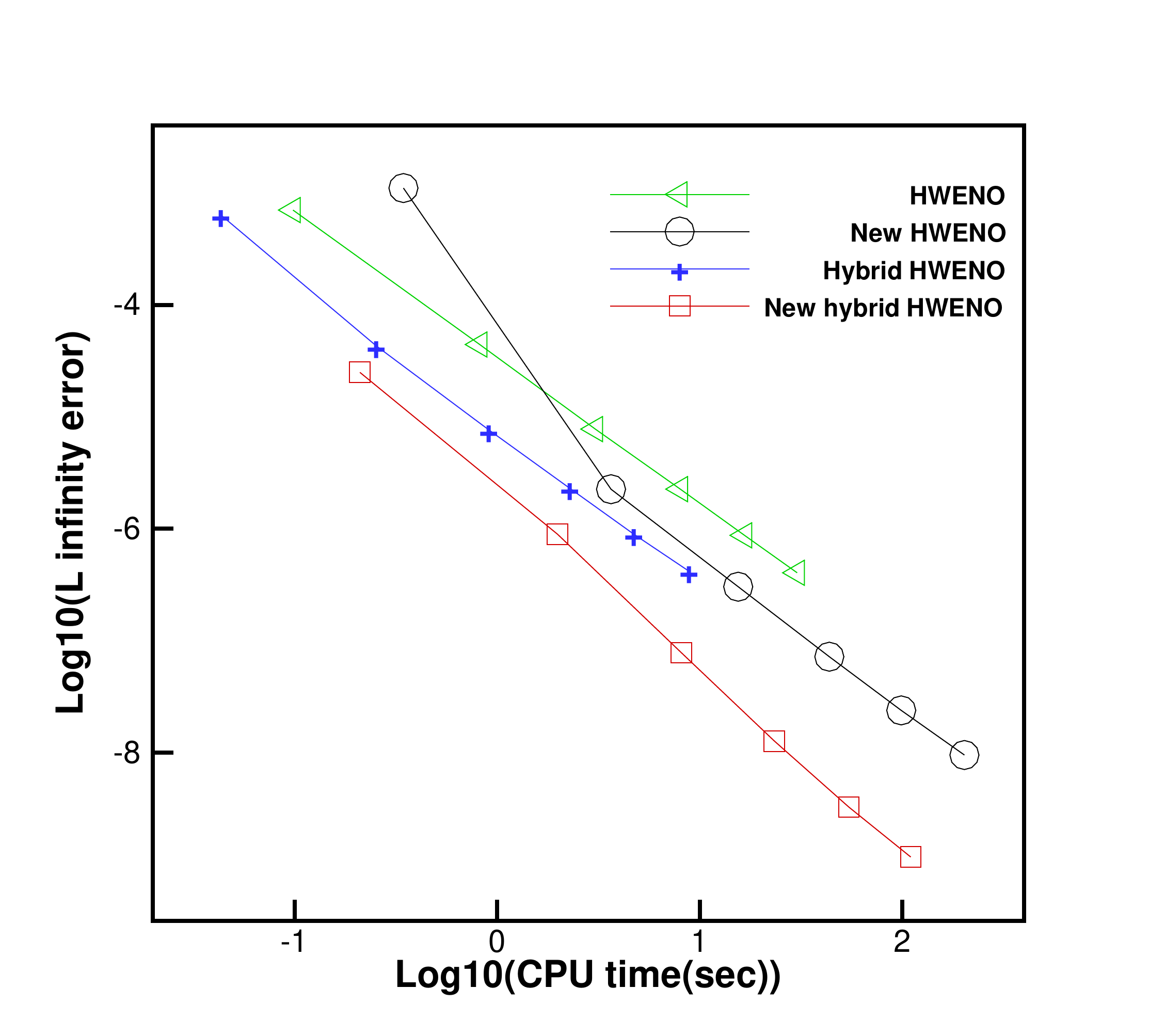,width=2.5 in}}
 \caption{2D-Burgers' equation: initial data
$u(x,y,0)=0.5+sin(\pi (x+y)/2)$. $T=0.5/\pi$. Computing times and errors. Triangle signs and a green solid  line: the results of HWENO scheme; circle signs and a black solid  line: the results of New HWENO scheme; plus signs and a blue solid  line: the results of Hybrid HWENO scheme; rectangle signs and a red solid  line: the results of New hybrid HWENO scheme.}
\label{Fburges2d_smooth}
\end{figure}

\noindent{\bf Example 3.3.} Two dimensional Burgers' equation:
\begin{equation}\label{2dbugers}
  u_t+(\frac {u^2} 2)_x+(\frac {u^2} 2)_y=0, \quad 0<x<4, \ 0<y<4.
\end{equation}
The initial condition is $u(x,y,0)=0.5+sin(\pi (x+y)/2)$ and periodic boundary conditions are applied in each direction. We compute the solution up to $T=0.5/\pi$, where the solution is smooth, and we present the numerical errors and orders in Table \ref{tburgers2d}, which illustrates that  New HWENO and New hybrid HWENO  schemes have the fifth order accuracy, while the  HWENO and hybrid HWENO schemes only have the fourth order accuracy, and we can see that different choice of the linear weights has no influence on the numerical accuracy for the new HWENO methodology. In addition, we present the numerical errors against CPU times by these HWENO schemes in Figure \ref{Fburges2d_smooth}, which illustrates New hybrid HWENO scheme has higher efficiency than Hybrid HWENO scheme with smaller numerical errors and higher order numerical accuracy, and the hybrid schemes both have less CPU times than the corresponding schemes. Meanwhile, New HWENO scheme has higher efficiency than HWENO scheme.
\begin{table}
\begin{center}
\caption{2D-Burgers' equation: initial data
$u(x,y,0)=0.5+sin(\pi (x+y)/2)$. HWENO  schemes. $T=0.5/\pi$. $L^1$ and $L^\infty$ errors and orders. Uniform meshes with $N_x\times N_y$ cells.}
\medskip
\begin{tabular} {lllllllll}
\hline
  $N_x\times N_y$ cells & \multicolumn{4}{l}{HWENO  scheme} & \multicolumn{4}{l}{New HWENO scheme}\\
\cline{2-5}\cline{6-9}
&$ L^1$ error &  order & $L^\infty$error &  order &$ L^1$ error &  order & $L^\infty$ error &order\\   \hline
$     40\times     40$ &     8.21E-05 &        &     7.02E-04 &       									&     1.28E-04 &        &     1.10E-03 &       \\
$     80\times     80$ &     4.67E-06 &       4.14 &     4.42E-05 &       3.99 									&     2.86E-07 &       8.81 &     2.25E-06 &       8.94 \\
$    120\times    120$ &     8.70E-07 &       4.15 &     7.76E-06 &       4.29 									&     2.52E-08 &       5.99 &     3.04E-07 &       4.95 \\
$    160\times    160$ &     2.66E-07 &       4.13 &     2.26E-06 &       4.29 									&     5.60E-09 &       5.22 &     7.19E-08 &       5.00 \\
$    200\times    200$ &     1.06E-07 &       4.12 &     8.73E-07 &       4.26 									&     1.79E-09 &       5.12 &     2.39E-08 &       4.95 \\
$    240\times    240$ &     5.02E-08 &       4.09 &     4.04E-07 &       4.23 									&     7.12E-10 &       5.05 &     9.53E-09 &       5.03 \\ \hline
  $N_x\times N_y$ cells & \multicolumn{4}{l}{Hybrid HWENO scheme} & \multicolumn{4}{l}{New Hybrid HWENO scheme}\\
\cline{2-5}\cline{6-9}
&$ L^1$ error &  order & $L^\infty$error &  order &$ L^1$ error &  order & $L^\infty$ error &order\\   \hline
$     40\times     40$ &     7.03E-05 &      &     6.32E-04 & 	&     2.70E-06 &   &  2.49E-05 &      \\
$     80\times     80$ &     3.93E-06 &       4.16 &     4.28E-05 &       3.88 									&     5.01E-08 &       5.75 &     8.91E-07 &       4.81 \\
$    120\times    120$ &     7.27E-07 &       4.16 &     7.61E-06 &       4.26 									&     4.15E-09 &       6.14 &     7.81E-08 &       6.00 \\
$    160\times    160$ &     2.18E-07 &       4.19 &     2.30E-06 &       4.16 									&     7.00E-10 &       6.18 &     1.27E-08 &       6.32 \\
$    200\times    200$ &     8.61E-08 &       4.16 &     8.95E-07 &       4.23 									&     1.94E-10 &       5.74 &     3.26E-09 &       6.09 \\
$    240\times    240$ &     4.05E-08 &       4.14 &     4.18E-07 &       4.18 									&     7.65E-11 &       5.12 &     1.17E-09 &       5.63 \\
\hline
\end{tabular}
\label{tburgers2d}
\end{center}
\end{table}
\smallskip

\begin{figure}
 \centerline{
\psfig{file=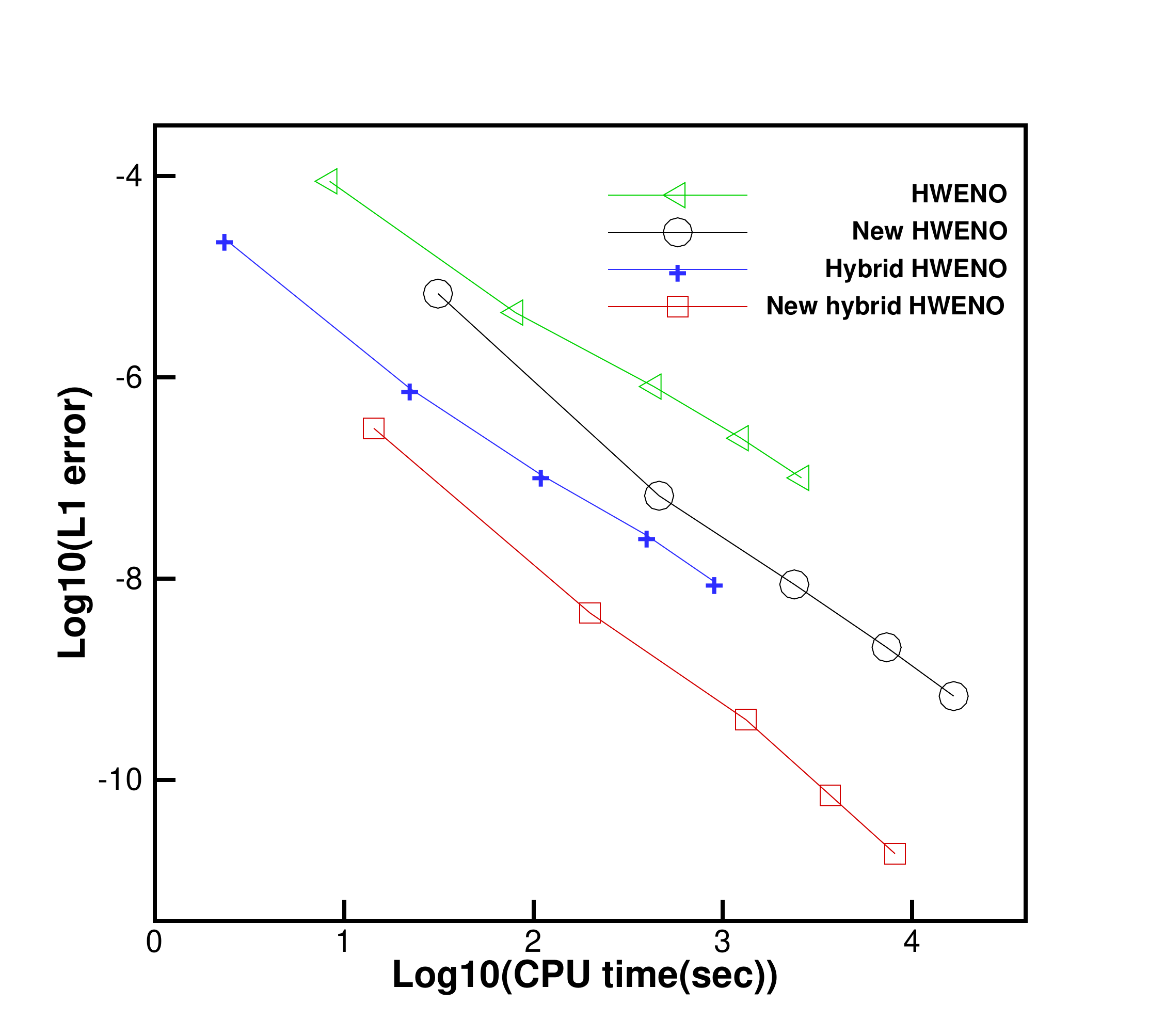,width=2.5 in}\psfig{file=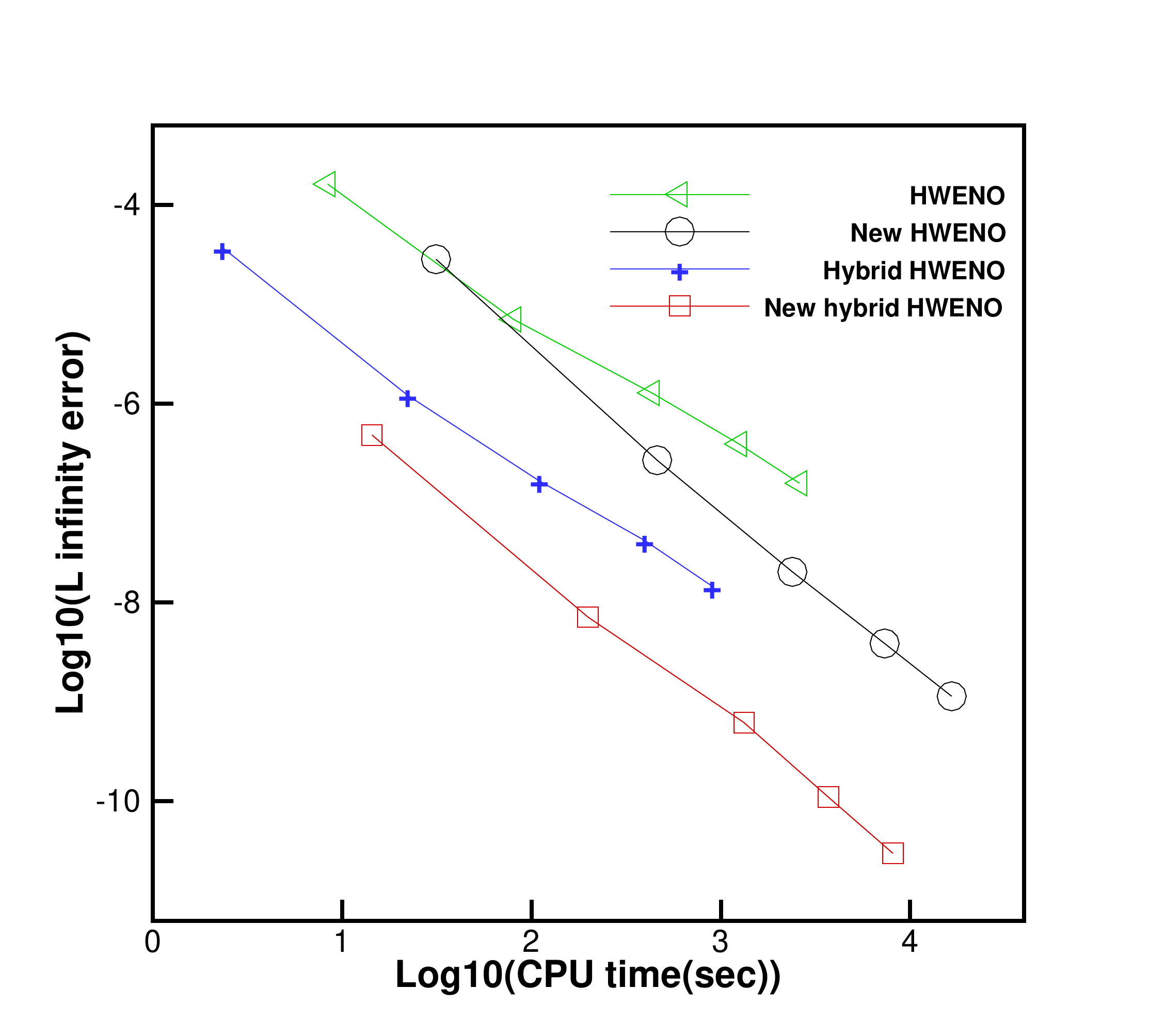,width=2.5 in}}
 \caption{2D-Euler equations: initial data
$\rho(x,y,0)=1+0.2\sin(\pi(x+y))$, $\mu(x,y,0)=1$,
$\nu(x,y,0)=1$ and $p(x,y,0)=1$. $T=2$. Computing times and errors. Triangle signs and a green solid  line: the results of HWENO scheme; circle signs and a black solid  line: the results of New HWENO scheme; plus signs and a blue solid  line: the results of Hybrid HWENO scheme; rectangle signs and a red solid  line: the results of New hybrid HWENO scheme.}
\label{FEuler2d}
\end{figure}

\noindent {\bf Example 3.4.} Two dimensional Euler equations:
\begin{equation}
\label{euler2}
 \frac{\partial}{\partial t}
 \left(
 \begin{array}{c}
 \rho \\
 \rho \mu \\
 \rho \nu \\
 E
 \end{array} \right )
 +
 \frac{\partial}{\partial x}
 \left (
\begin{array}{c}
\rho \mu \\
\rho \mu^{2}+p \\
\rho \mu \nu \\
\mu(E+p)
\end{array}
 \right )
 +
 \frac{\partial}{\partial y}
 \left (
\begin{array}{c}
\rho \nu \\
\rho \mu \nu \\
\rho \nu^{2}+p \\
\nu(E+p)
\end{array}
 \right )=0,
\end{equation}
in which $\rho$  is the density; $(\mu,\nu)$ is the velocity; $E$ is the total energy; and $p$ the is pressure. The initial conditions are $\rho(x,y,0)=1+0.2\sin(\pi(x+y))$, $\mu(x,y,0)=1$, $\nu(x,y,0)=1$, $p(x,y,0)=1$ and $\gamma=1.4$. The computing domain is $(x,y)\in [0,2] \times [0, 2]$ with periodic boundary conditions in $x$ and $y$ directions, respectively. The exact solution of $\rho$ is $\rho(x,y,t)=1+0.2\sin(\pi(x+y-2t))$ and the computing time is  $T=2$. We give the numerical errors and  orders of the density for HWENO, New HWENO, Hybrid HWENO, New hybrid HWENO schemes in Table \ref{tEluer2d}, then, we can find the New HWENO and New hybrid HWENO achieve the fifth order accuracy, but the HWENO and Hybrid HWENO scheme only have the fourth order accuracy, meanwhile, we can see that  random positive linear weights (the sum equals one) would have no impact on the order accuracy of New HWENO scheme. Finally, we also show their numerical errors against CPU times in Figure \ref{FEuler2d}, which illustrates  New hybrid HWENO scheme has  higher efficiency than other three schemes, meanwhile, New  HWENO scheme has better performance with less numerical errors and higher order accuracy than HWENO scheme.
\begin{table}
\begin{center}
\caption{2D-Euler equations: initial data $\rho(x,y,0)=1+0.2\sin(\pi(x+y))$, $\mu(x,y,0)=1$,
$\nu(x,y,0)=1$ and $p(x,y,0)=1$. HWENO schemes. $T=2$. $L^1$ and $L^\infty$ errors and orders. Uniform meshes with $N_x\times N_y$ cells.}
\medskip
\begin{tabular} {lllllllll}
\hline
$N_x\times N_y$ cells & \multicolumn{4}{l}{HWENO  scheme} & \multicolumn{4}{l}{New HWENO scheme}\\
\cline{2-5}\cline{6-9}
  &$ L^1$ error &  order & $L^\infty$error &  order &$ L^1$ error &  order & $L^\infty$ error &order\\   \hline
$     30\times     30$ &     8.85E-05 &        &     1.63E-04 &       							 &     6.78E-06 &       &     2.86E-05 &       \\
$     60\times     60$ &     4.39E-06 &       4.33 &     7.09E-06 &       4.52 									 &     6.64E-08 &       6.67 &     2.71E-07 &       6.72 \\
$     90\times     90$ &     8.08E-07 &       4.17 &     1.29E-06 &       4.20 									 &     8.73E-09 &       5.00 &     2.04E-08 &       6.38 \\
$    120\times    120$ &     2.48E-07 &       4.11 &     3.95E-07 &       4.11 									 &     2.07E-09 &       5.00 &     3.88E-09 &       5.77 \\
$    150\times    150$ &     1.00E-07 &       4.07 &     1.59E-07 &       4.07 									 &     6.78E-10 &       5.00 &     1.14E-09 &       5.48 \\
\hline
$N_x\times N_y$ cells & \multicolumn{4}{l}{Hybrid HWENO scheme} & \multicolumn{4}{l}{New hybrid HWENO scheme}\\
\cline{2-5}\cline{6-9}
&$ L^1$ error &  order & $L^\infty$error &  order &$ L^1$ error &  order & $L^\infty$ error &order\\   \hline
$     30\times     30$ &     2.37E-05 &        &     3.72E-05 &       									&     3.11E-07 &       &     4.87E-07 &        \\
$     60\times     60$ &     7.76E-07 &       4.93 &     1.22E-06 &       4.93									&     4.55E-09 &       6.09 &     7.14E-09 &       6.09 \\
$     90\times     90$ &     1.07E-07 &       4.89 &     1.68E-07 &       4.89									&     3.95E-10 &       6.03 &     6.20E-10 &       6.03 \\
$    120\times    120$ &     2.67E-08 &       4.82 &     4.20E-08 &       4.82									&     7.01E-11 &       6.01 &     1.10E-10 &       6.00 \\
$    150\times    150$ &     9.27E-09 &       4.75 &     1.46E-08 &       4.75									&     1.84E-11 &       5.99 &     3.00E-11 &       5.84 \\
\hline
\end{tabular}
\label{tEluer2d}
\end{center}
\end{table}
\smallskip

\subsection{Non-smooth tests}

We present the results of the hybrid HWENO scheme here, meanwhile, the linear weights for the low degree polynomials are set as 0.01 and the linear weight for the high degree polynomial is the rest (the sum of their linear weights equals one). For comparison, we also show the numerical results of the hybrid HWENO scheme \cite{ZCQHH}. From the results of the non-smooth tests, two schemes have similar performances in one dimension, but the new hybrid HWENO scheme has better numerical performances in two dimension for the new hybrid HWENO scheme has higher order numerical accuracy. In addition, the new hybrid HWENO scheme uses more simpler HWENO methodology, where any artificial positive linear weights (the sum equals 1) can be used, which is easier to implement in the computation, and it also uses less candidate stencils and bigger CFL number for two dimensional problems.

\begin{figure}
 \centerline{ \psfig{file=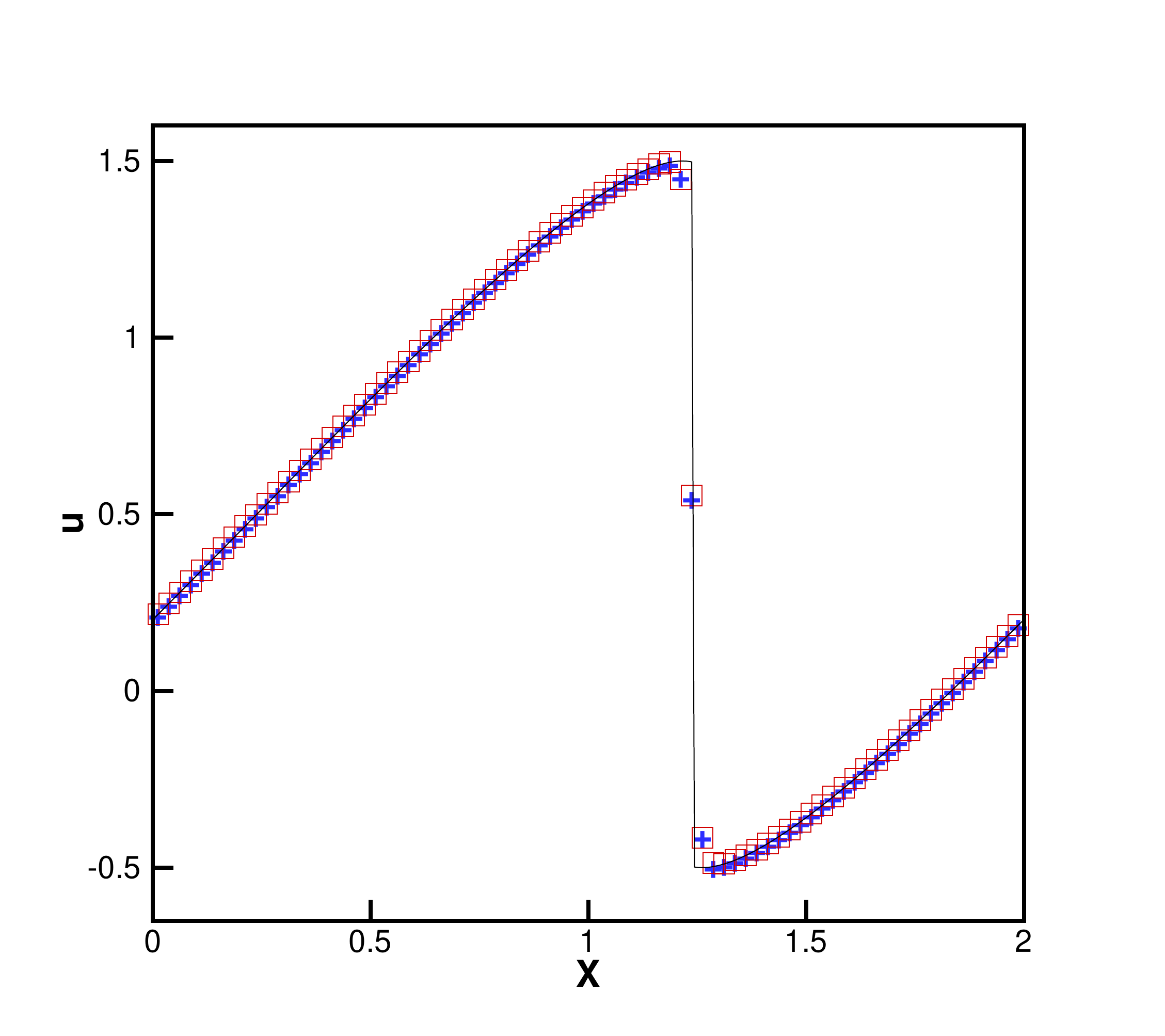,width=3 in}}
 \caption{1D-Burgers' equation: initial data
$u(x,0)=0.5+sin(\pi x)$. $T=1.5/\pi$. Black solid line: exact solution; blue plus signs: the results of the hybrid HWENO scheme; red squares: the results of the new hybrid HWENO scheme. Uniform meshes with 80 cells.}
\label{Fburges1d}
\end{figure}

\noindent{\bf Example 3.5.} We solve the one-dimensional Burgers' equation (\ref{1dbugers}) as introduced in Example 3.1 with  same initial and boundary conditions, but the final computing time is $t=1.5/\pi$, in which the solution is discontinuous. In Figure \ref{Fburges1d}, we present the the numerical solution of the HWENO schemes and the exact solution, and we can see that two schemes have similar numerical results with high resolutions.
\smallskip

\begin{figure}
\centering{\psfig{file=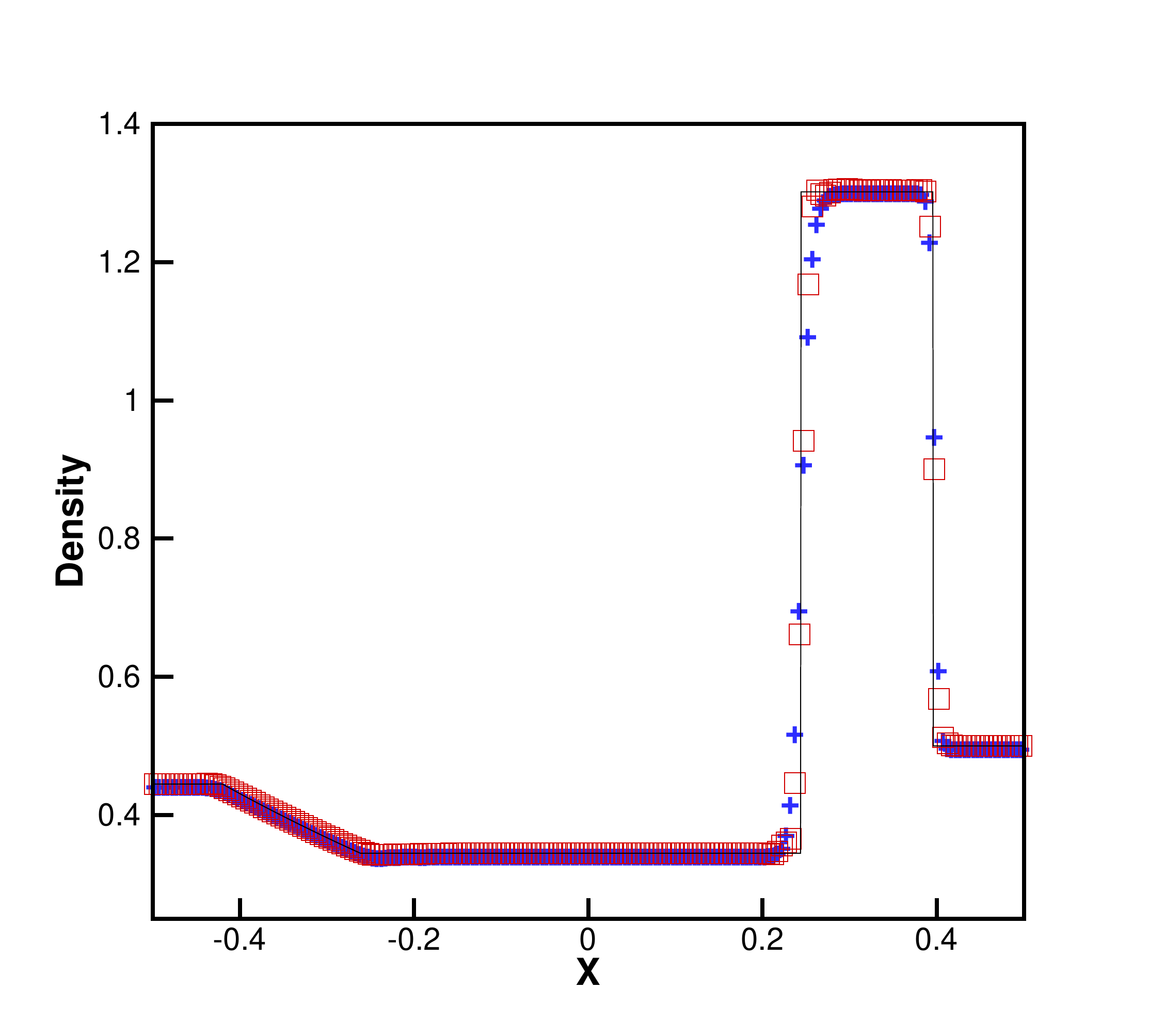,width=2 in}\psfig{file=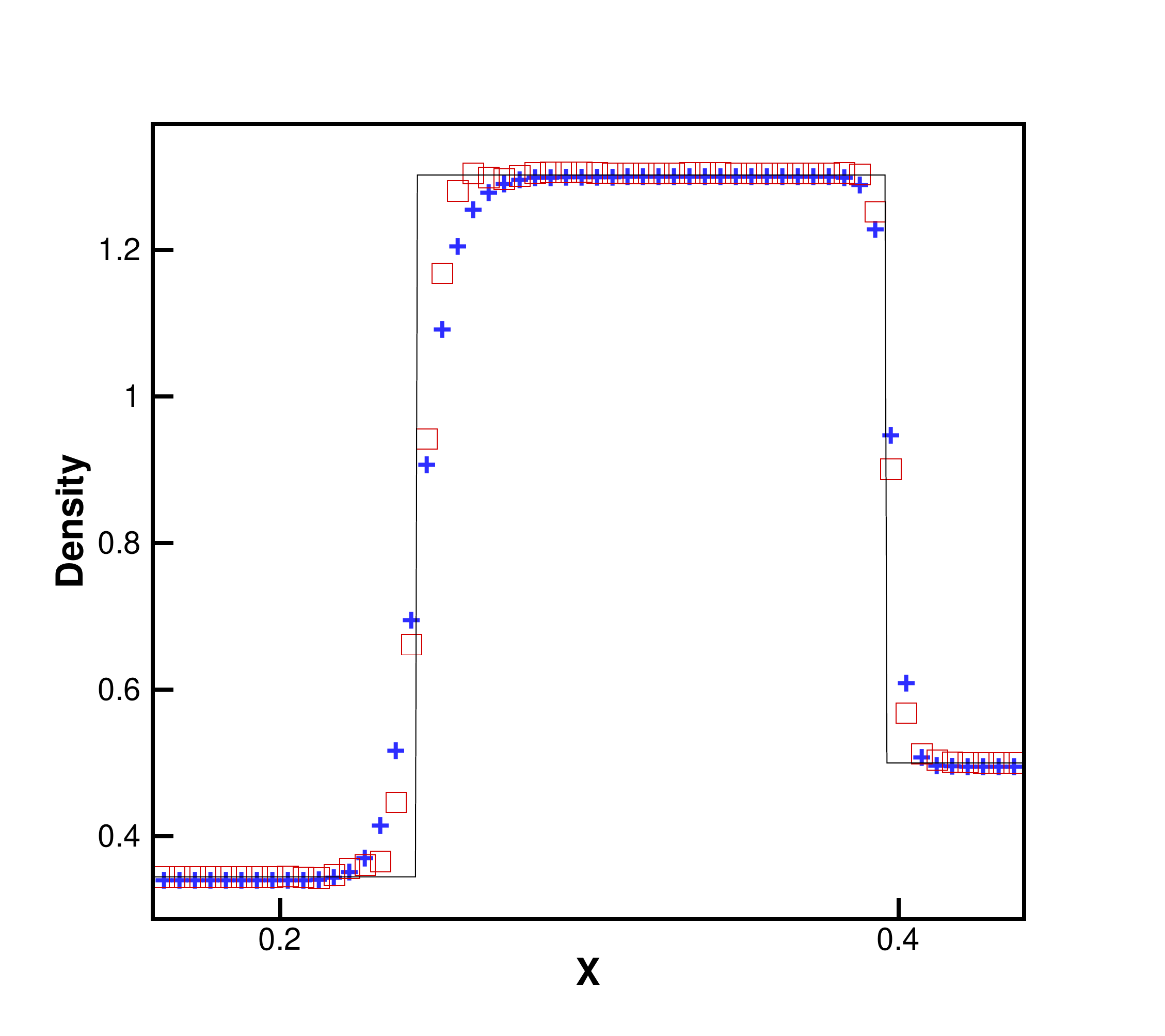,width=2 in}\psfig{file=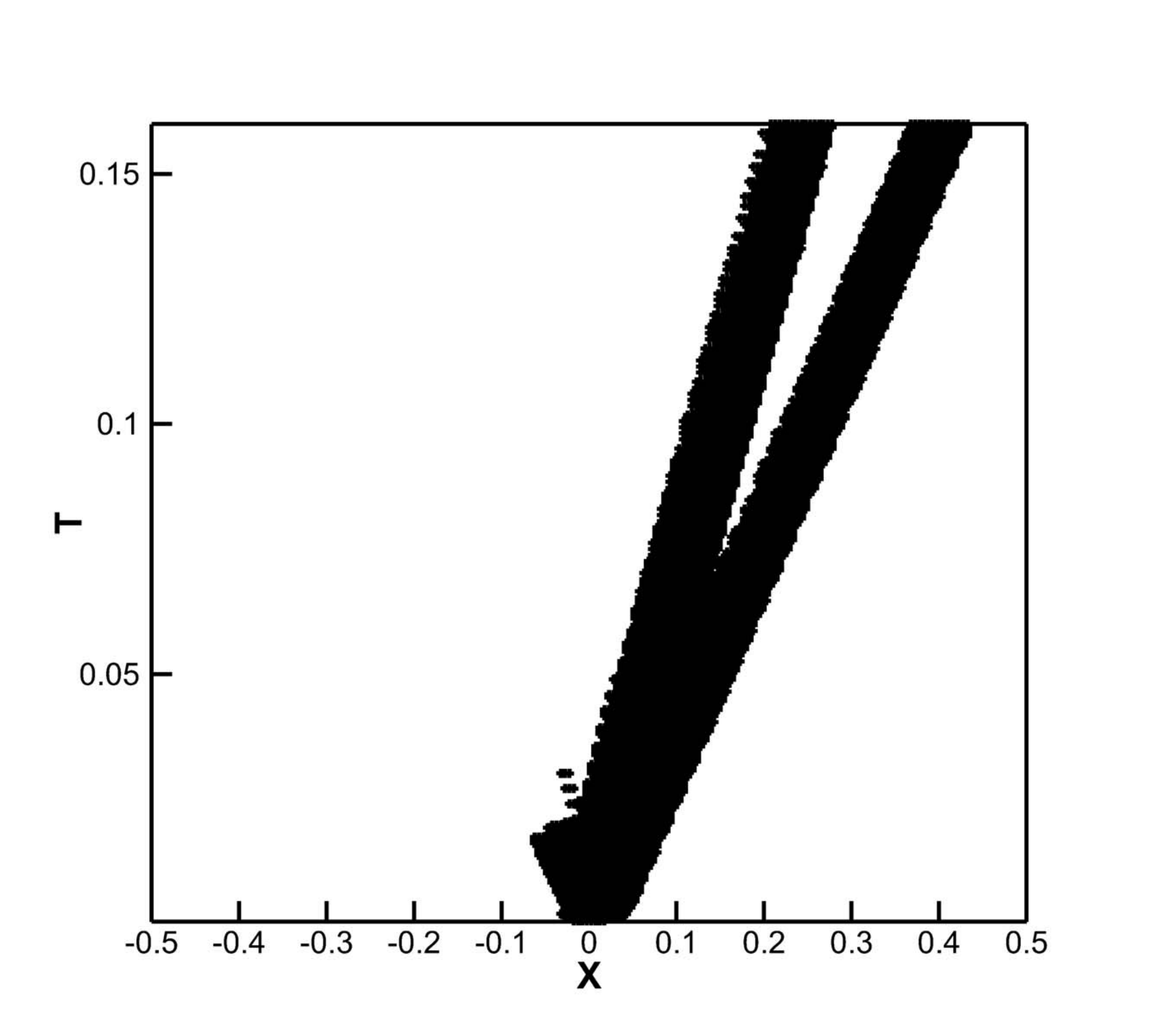,width=2 in}}
 \caption{The Lax problem.  T=0.16. From left to right: density; density zoomed in; the cells where  the modification for the first order moments are computed in the new hybrid HWENO scheme. Black solid line: the exact solution; blue plus signs: the results of the hybrid HWENO scheme; red squares: the results of the new hybrid HWENO scheme. Uniform meshes with 200 cells.}
\label{laxfig}
\end{figure}

\noindent {\bf Example 3.6.} The Lax problem for 1D Euler equations with the next Riemann initial condition:
\begin{equation*}
\label{lax} (\rho,\mu,p,\gamma)^T= \left\{
\begin{array}{ll}
(0.445,0.698,3.528,1.4)^T,& x \in [-0.5,0),\\
(0.5,0,0.571,1.4)^T,& x \in [0,0.5].
\end{array}
\right.
\end{equation*}
The computing time is $T=0.16$. In Figure \ref{laxfig}, we plot  the exact solution against the computed density $\rho$ obtained with the HWENO schemes, the zoomed in picture and the time history of the cells where the modification procedure is used in the new hybrid  HWENO scheme. We can see the results computed by the new hybrid HWENO schemes is closer to the exact solution, and we also find that only 13.41 \% cells where we use the new HWENO methodology, which means that most regions directly use linear approximation with no modification for the first order moments and no HWENO reconstruction for the spatial discretization. The new hybrid HWENO scheme keeps good resolutions too.
\smallskip

\begin{figure}
 \centering{\psfig{file=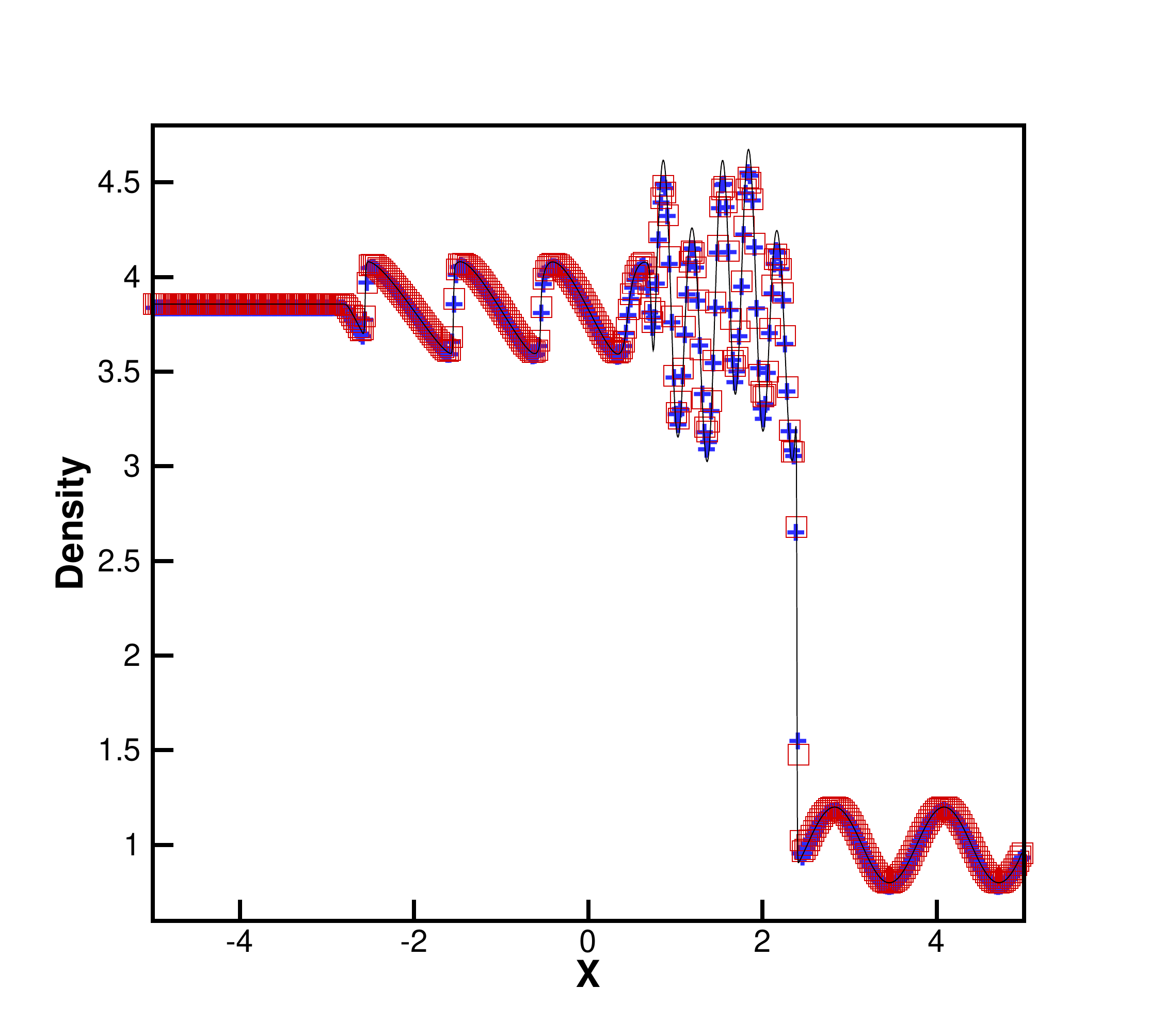,width=2 in}\psfig{file=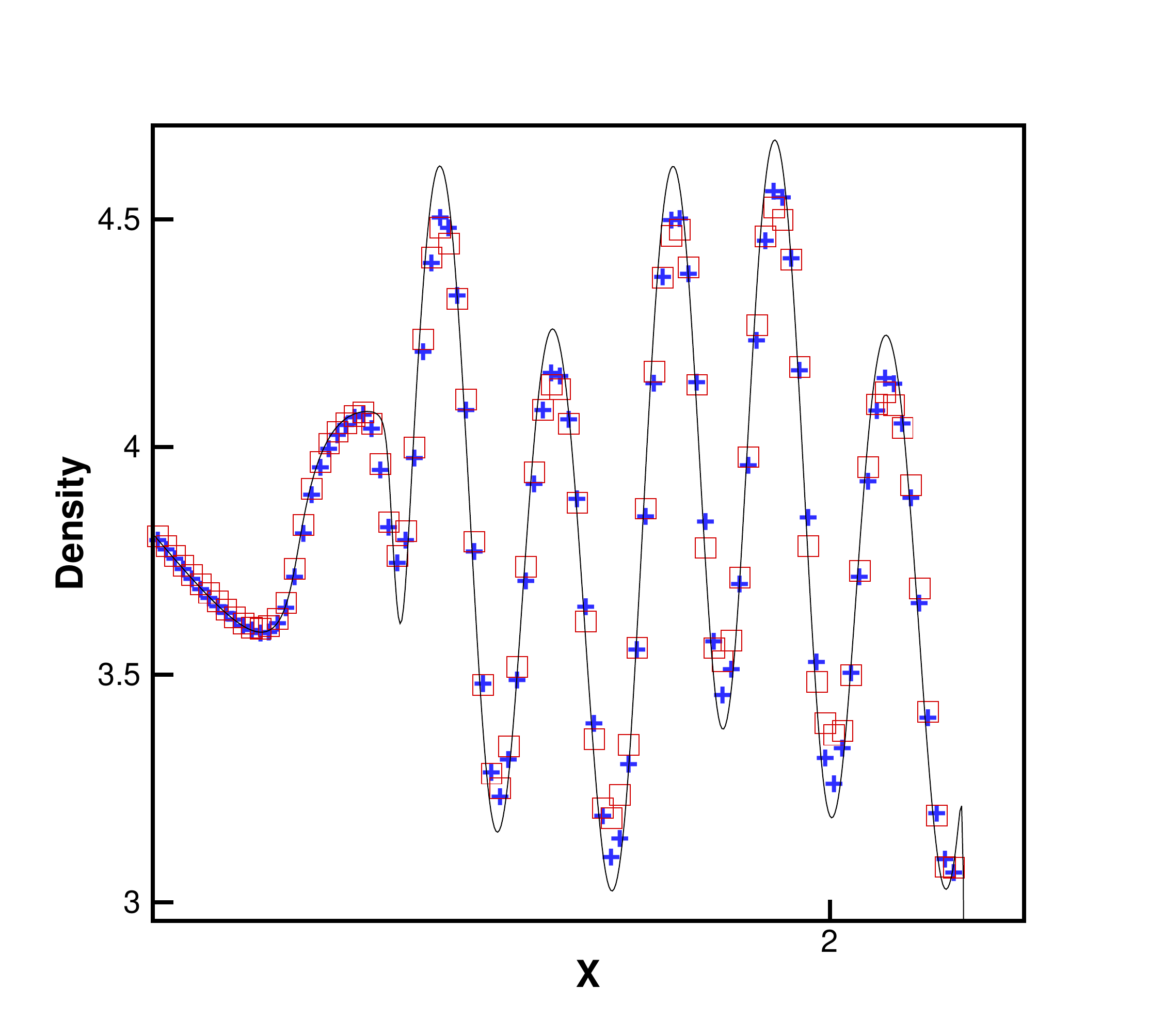,width=2 in}\psfig{file=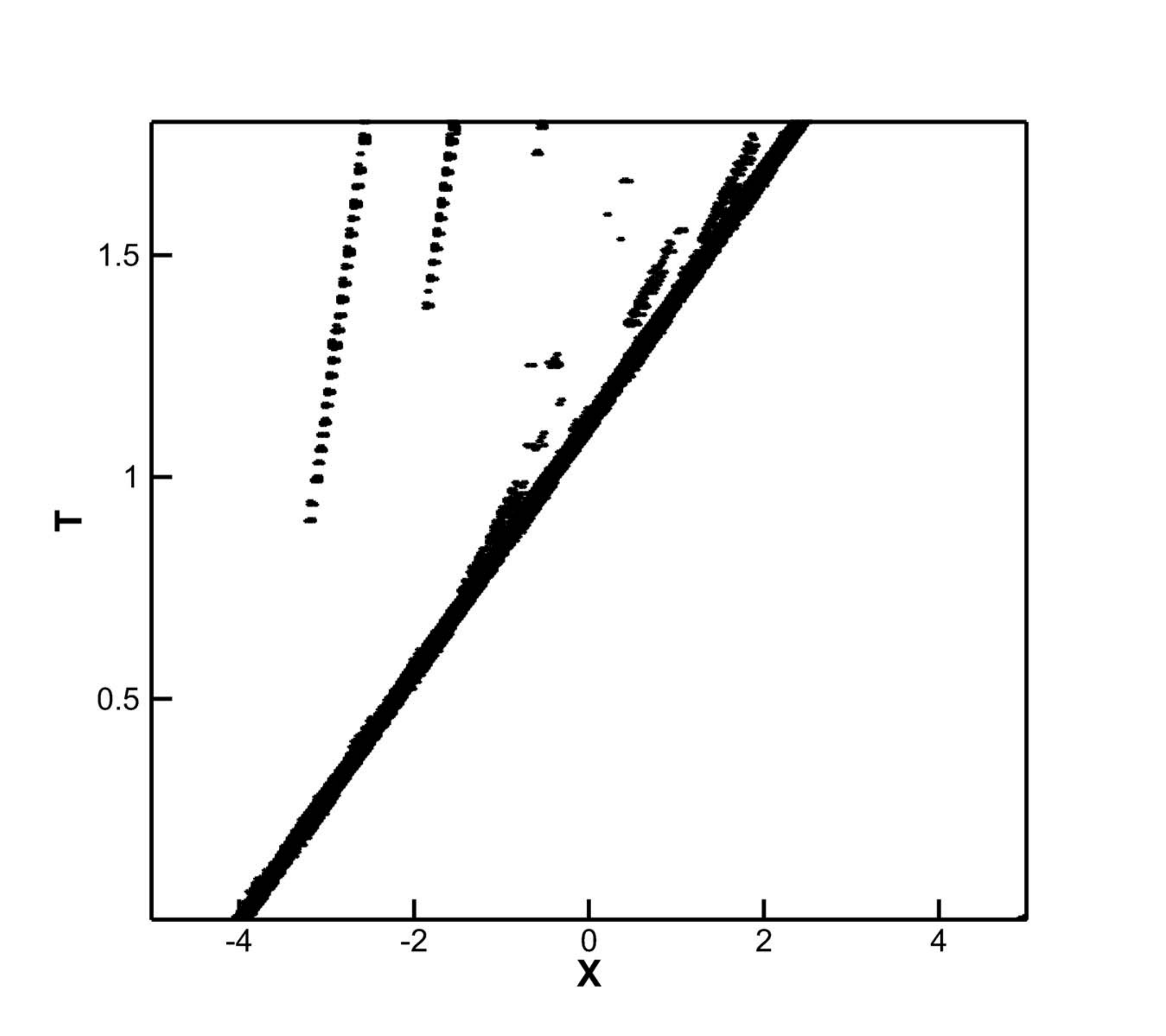,width=2 in}}
\caption{The shock density wave interaction problem. T=1.8. From left to right: density; density zoomed in; the cells where  the modification for the first order moments are computed in the new hybrid HWENO scheme. Black solid line: the exact solution; blue plus signs: the results of the hybrid HWENO scheme; red squares: the results of the new hybrid HWENO scheme. Uniform meshes with 400 cells.}
\label{sin}
\end{figure}

\noindent {\bf Example 3.7.} The Shu-Osher problem, which has a shock interaction with entropy waves \cite{s2}. The initial condition is
\begin{equation*}
\label{ShuOsher} (\rho,\mu,p,\gamma)^T= \left\{
\begin{array}{ll}
(3.857143, 2.629369, 10.333333,1.4)^T,& x \in [-5, -4),\\
(1 + 0.2\sin(5x), 0, 1,1.4)^T,& x \in [-4,5].
\end{array}
\right.
\end{equation*}
This is a typical example both containing  shocks and complex smooth region structures, which has a moving Mach=3 shock interacting with sine waves in density. The computing time is up to $T=1.8$. In Figure \ref{sin}, we plot the computed density $\rho$ by HWENO schemes against the referenced "exact" solution, the zoomed in picture and the time history of the troubled-cells for the new hybrid HWENO scheme. The referenced "exact" solution is computed by the fifth order finite difference WENO scheme \cite{js} with 2000 grid points. We can see two schemes have similar numerical results with high resolutions, but the new hybrid HWENO scheme doesn't need to calculate the linear weights in advance. In addition, only 3.54\% cells are identified as the troubled-cells where we need to modify their first order moments.
\smallskip

\begin{figure}
  \centering{\psfig{file=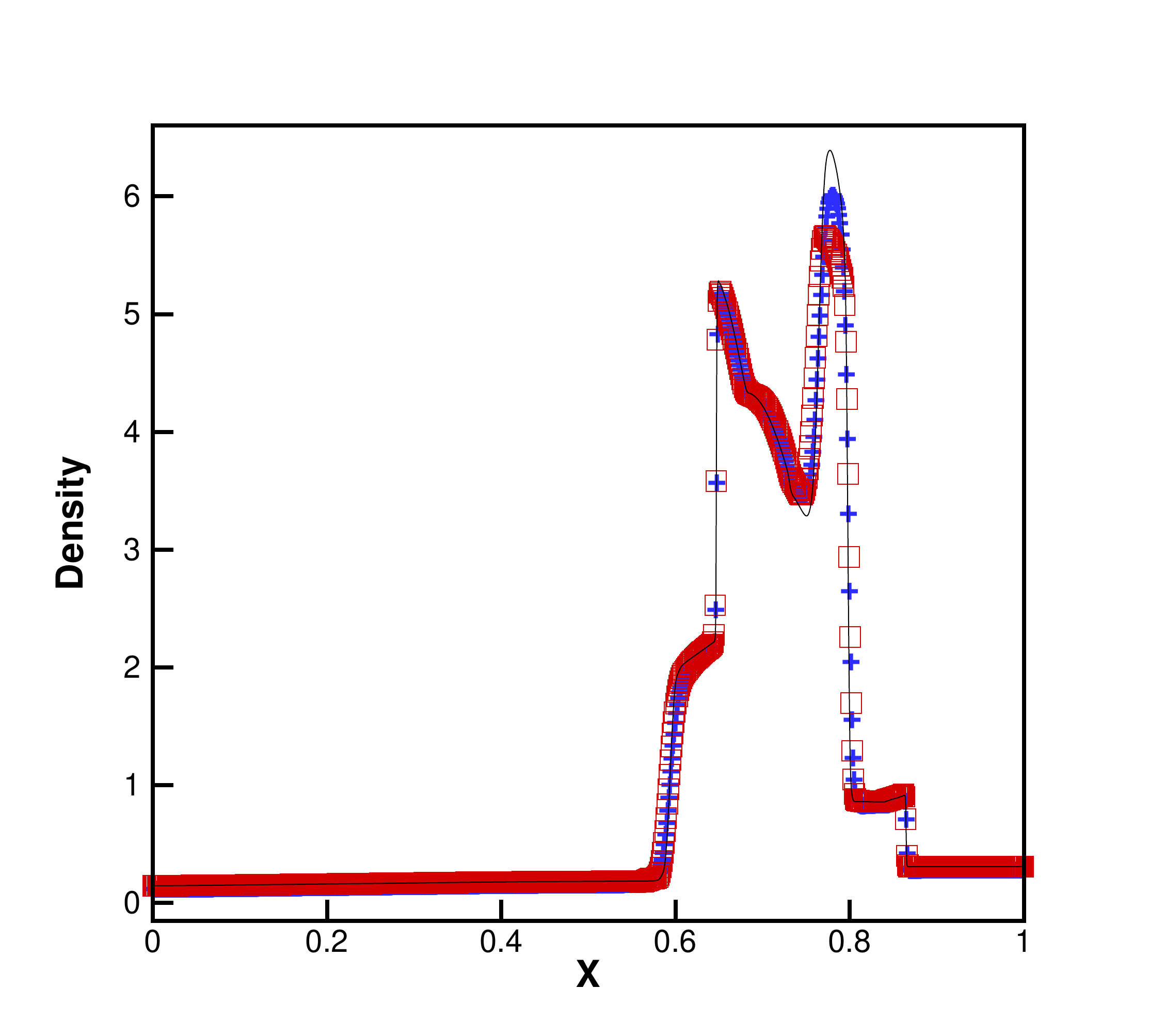,width=2 in}\psfig{file=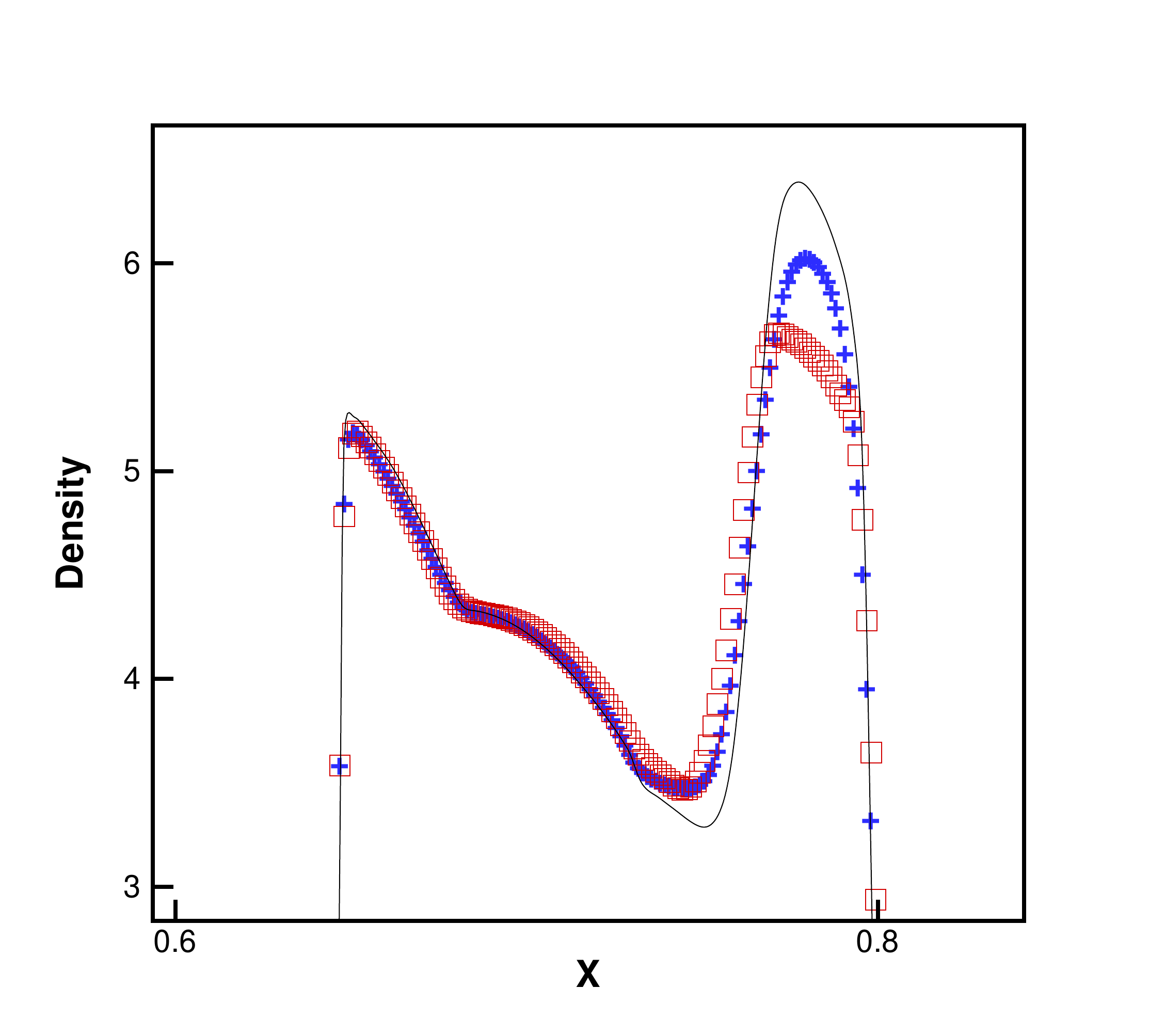,width=2 in}\psfig{file=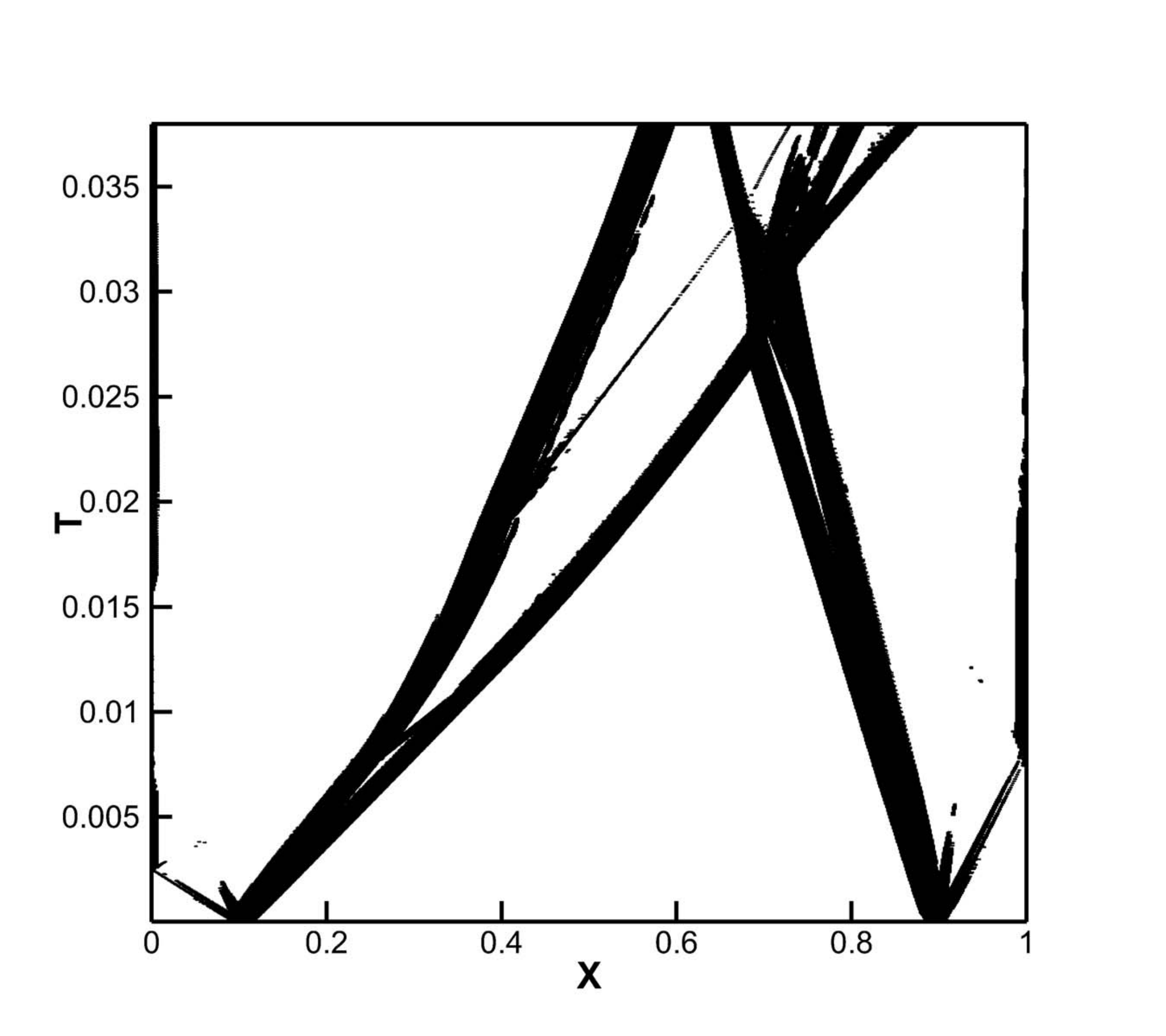,width=2 in}}
\caption{The blast wave problem. T=0.038. From left to right: density; density zoomed in; the cells where  the modification of the first order moments are computed in the new hybrid HWENO scheme. Black solid line: the exact solution; blue plus signs: the results of the hybrid HWENO scheme; red squares: the results of the new hybrid HWENO scheme. Uniform meshes with 800 cells.}
\label{blast}
\end{figure}

\noindent {\bf Example 3.8.} We solve the next interaction of two blast waves problems. The initial conditions are:
\begin{equation*}
\label{blastwave} (\rho,\mu,p,\gamma)^T= \left\{
\begin{array}{ll}
(1,0,10^3,1.4)^T,& 0<x<0.1,\\
(1,0,10^{-2},1.4)^T,& 0.1<x<0.9,\\
(1,0,10^2,1.4)^T,& 0.9<x<1.
\end{array}
\right.
\end{equation*}
The computing time is $T=0.038$, and the reflective boundary condition is applied here. In Figure \ref{blast}, we also plot the computed density against the reference "exact" solution, the zoomed in picture  and the time history of the troubled-cells. The reference "exact" solution is also computed by the fifth order finite difference WENO scheme \cite{js} with 2000 grid points. We notice that the hybrid HWENO scheme has better performance than the new hybrid HWENO scheme. The reason maybe that the modification for the first order moments uses more information provided by the two linear polynomials in this example, but the new HWENO methodology is easy to implement in the computation. Similarly, only 13.94\% cells are identified as the troubled-cells, and we directly use high order linear approximation on other cells.
\smallskip

\begin{figure}
 \centerline{\psfig{file=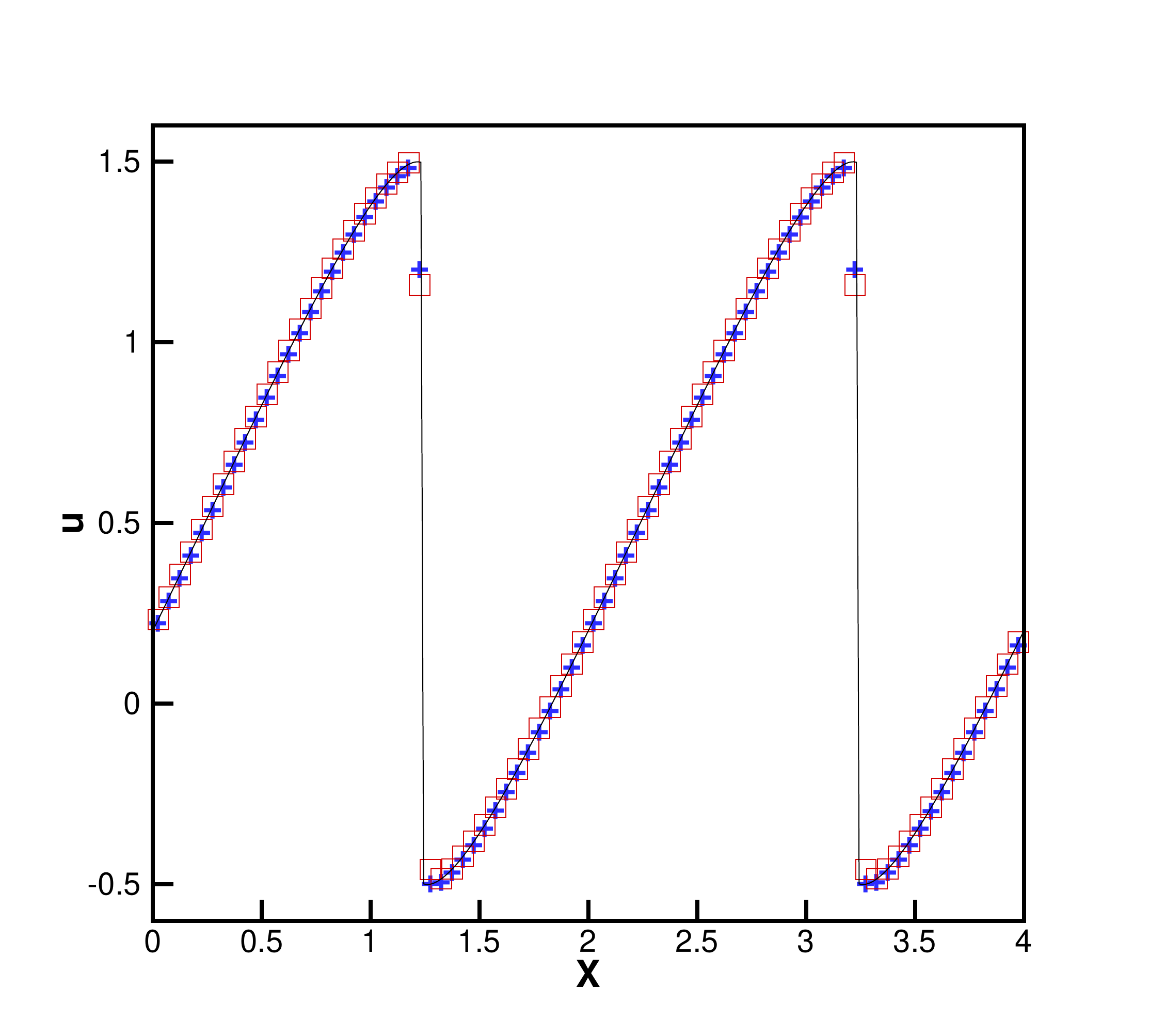,width=2.5 in}\psfig{file=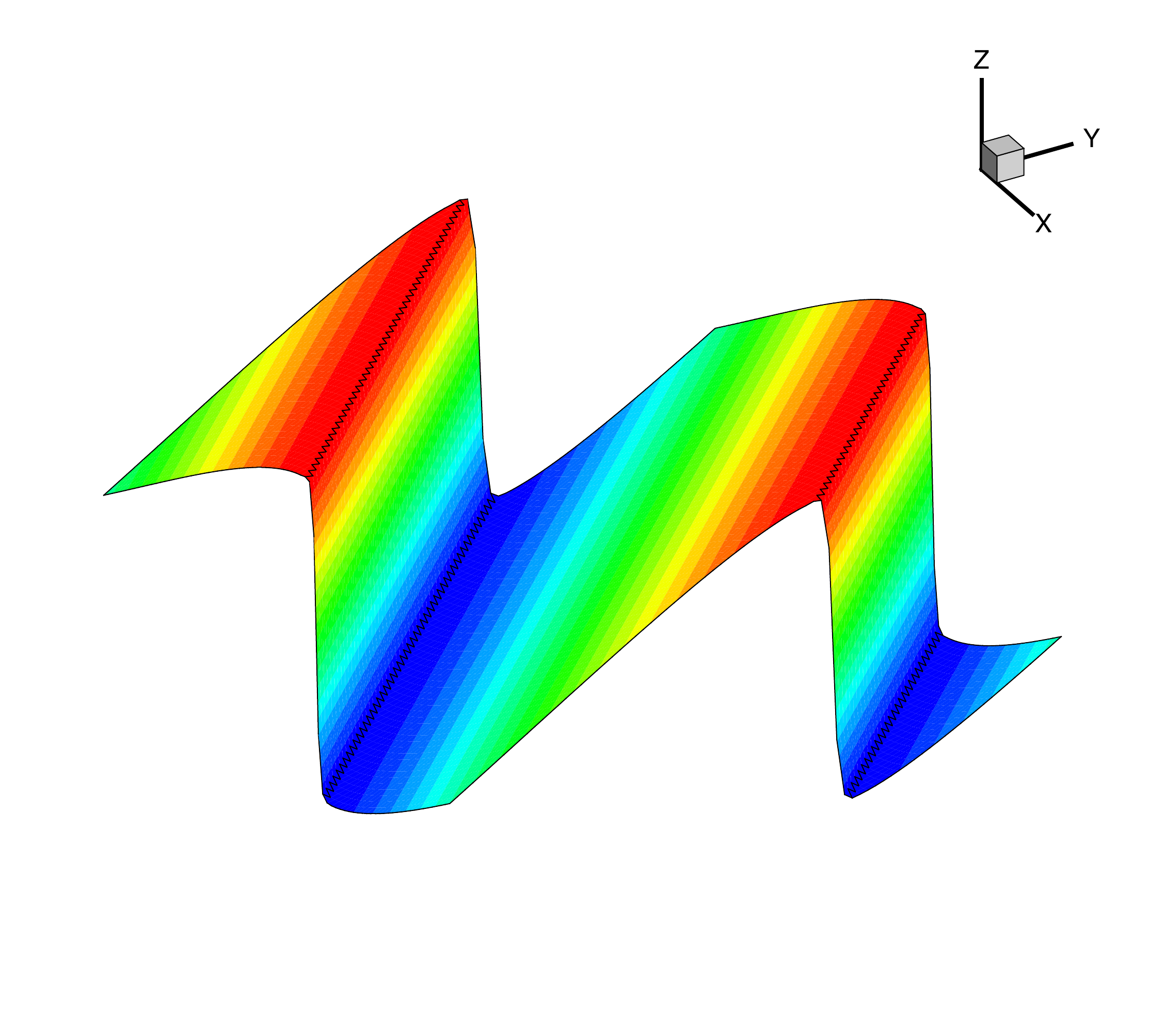,width=2.5 in}}
 \caption{2D-Burgers' equation: initial data
$u(x,y,0)=0.5+sin(\pi (x+y)/2)$. $T=1.5/\pi$. From left to right: the numerical solution at $x=y$ computed by HWENO schemes; the surface of the numerical solution for the new hybrid HWENO scheme. Black solid line: exact solution; blue plus signs: the results of the hybrid HWENO scheme; red squares: the results of the new hybrid HWENO scheme. Uniform meshes with $80\times80$ cells.}
\label{Fburges2d}
\end{figure}

\noindent{\bf Example 3.9.} We solve the two-dimensional Burgers' equation (\ref{2dbugers}) given in Example 3.3. The same initial and boundary conditions are applied here, but the computing time is up to $T=1.5/\pi$, in which the solution is discontinuous. In Figure \ref{Fburges2d}, we present the numerical solution computed by HWENO schemes against the exact solution and the surface of the numerical solution by the new hybrid HWENO scheme. Similarly, we can see the HWENO schemes have high resolutions.
\smallskip

\begin{figure}
\centerline{\psfig{file=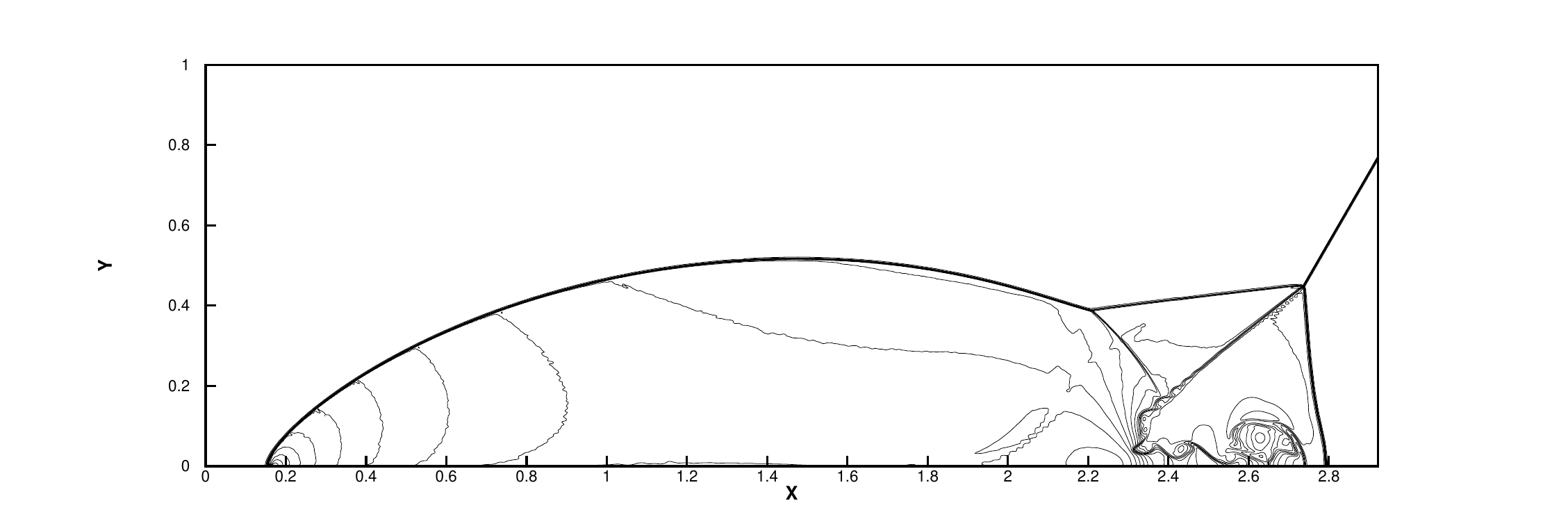,width=3.5 in} \psfig{file=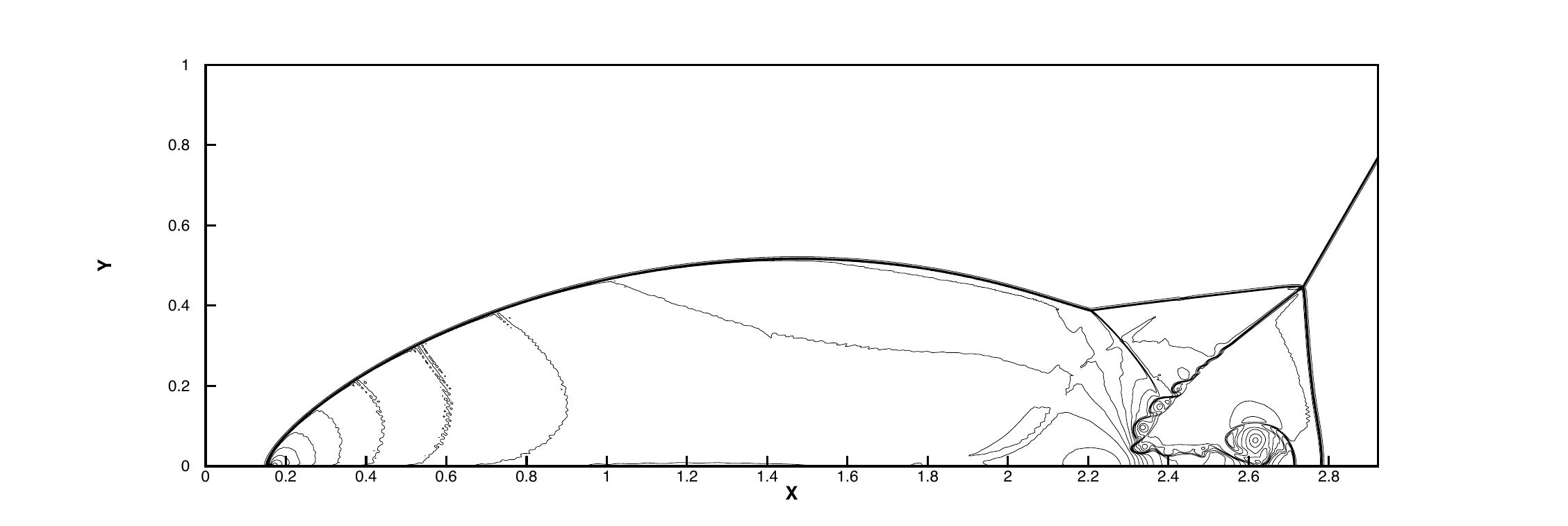,width=3.5 in}}
\centerline{\psfig{file=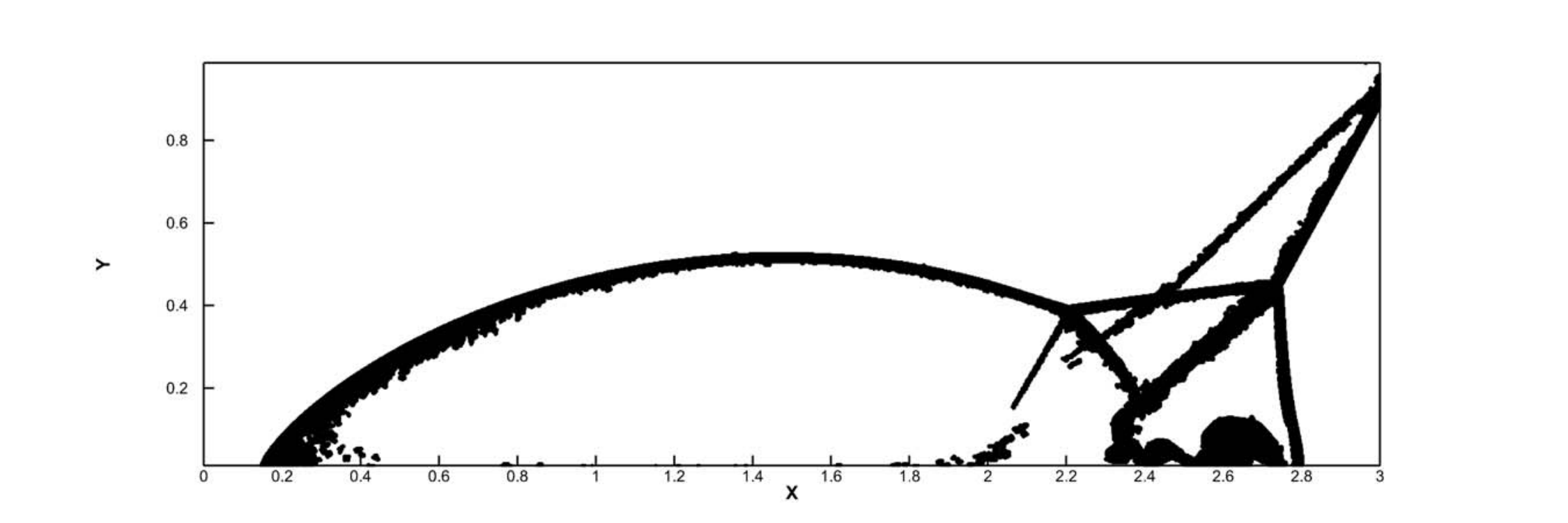,width=3.5 in} \psfig{file=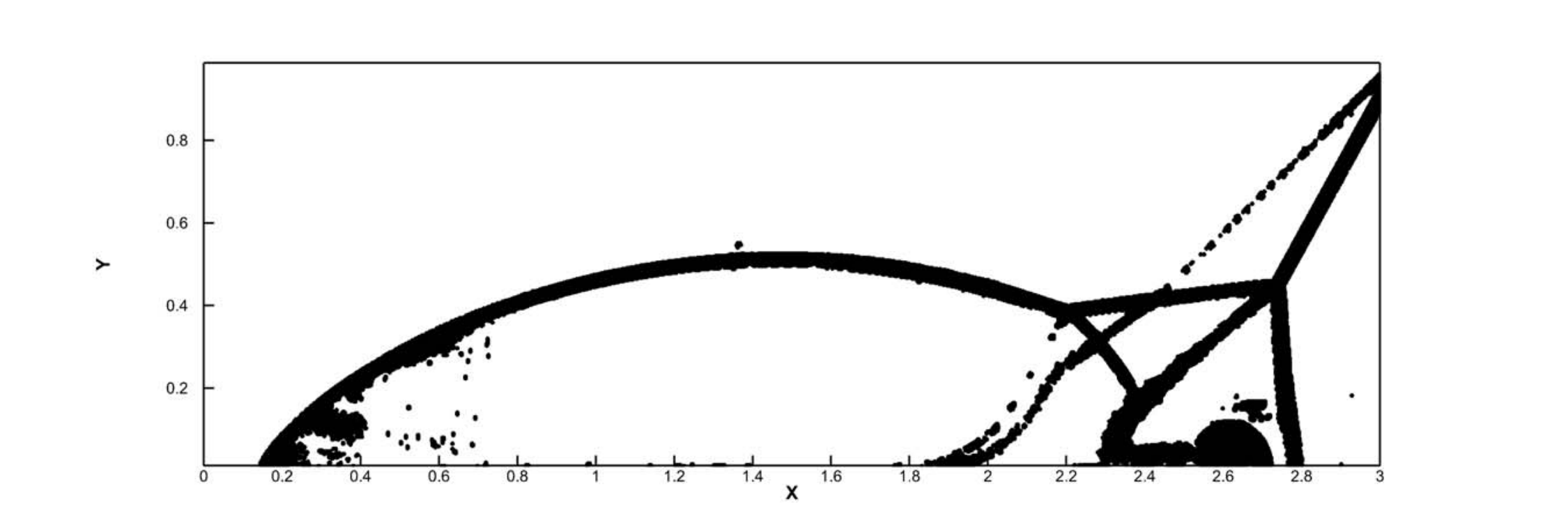,width=3.5 in}}
\centerline{\psfig{file=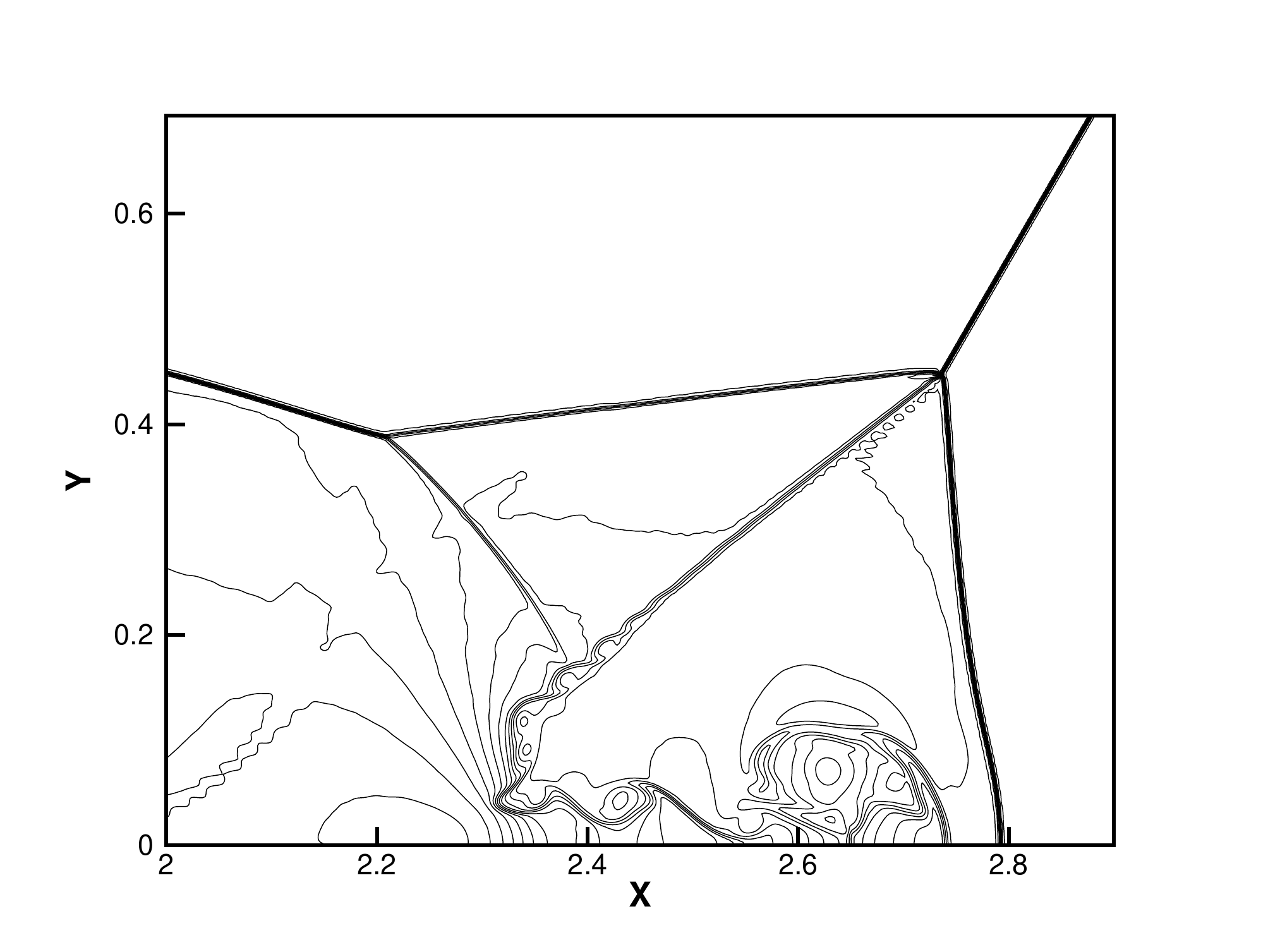,width=2 in} \psfig{file=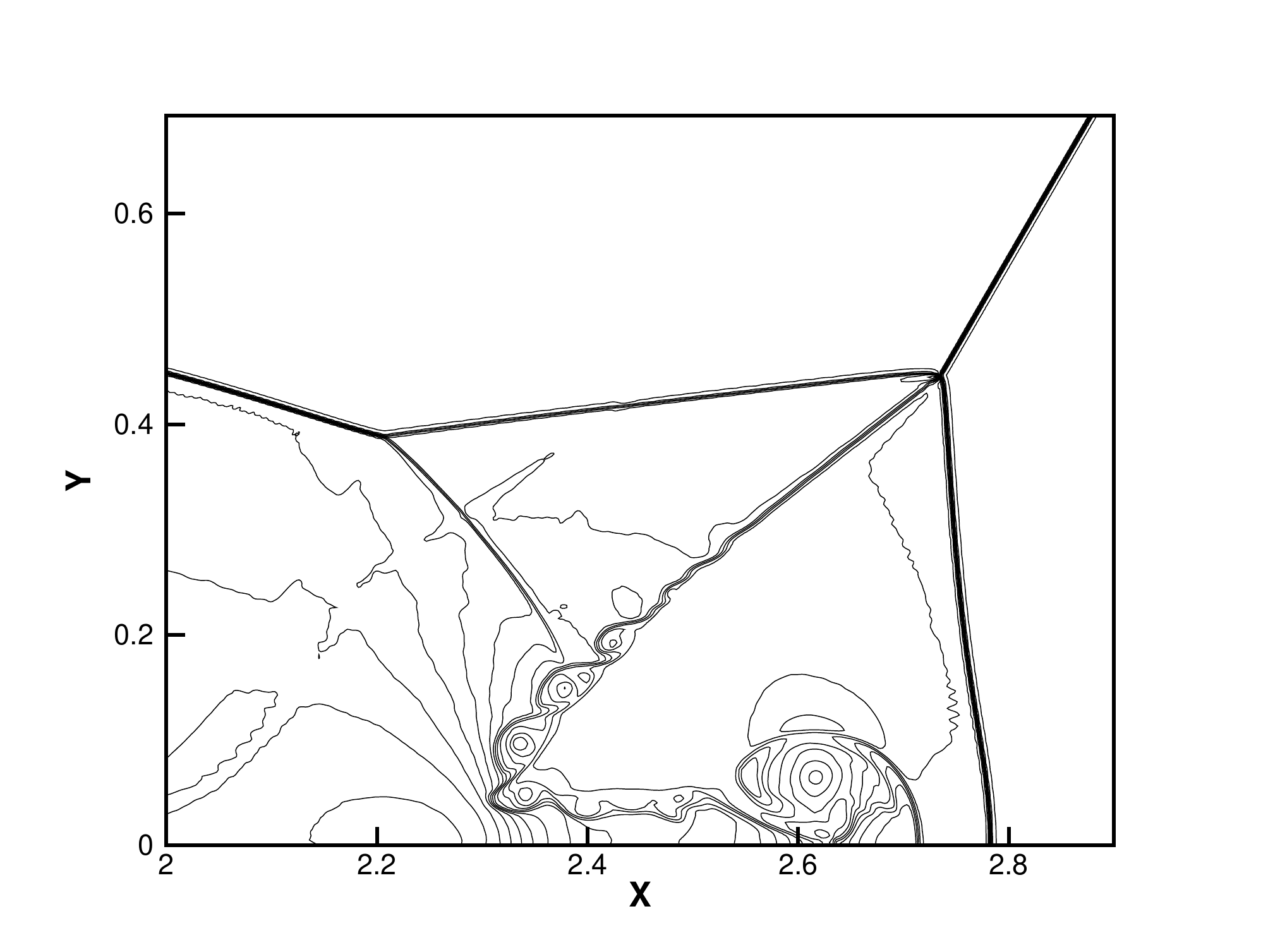,width=2 in} }
 \caption{ Double Mach reflection problem. T=0.2. From top to bottom: 30 equally spaced density contours from 1.5 to
 22.7; the locations of the troubled-cells at the final time; zoom-in pictures around the Mach stem. The hybrid HWENO scheme (left); the new hybrid HWENO scheme (right). Uniform meshes with 1920 $\times$ 480 cells.}
\label{smhfig}
\end{figure}

\noindent {\bf Example 3.10.} We now solve double Mach reflection problem \cite{Wooc} modeled by the two-dimensional Euler equations (\ref{euler2}). The computational domain is
\( [0,4]\times[0,1]\). The boundary conditions are: a reflection wall lies at the bottom from $x=\frac{1}{6}$, $y=0$ with a $60^{o}$ angle based on $x$-axis. For the bottom boundary, the reflection boundary condition are applied, but the part from $x=0$ to $x=\frac{1}{6}$ imposes the exact post-shock  condition. For the top boundary, it is  the exact motion of the Mach 10 shock.  $\gamma=1.4$ and the final computing time is up to $T=0.2$.  In Figure \ref{smhfig}, we plot the pictures of region \( [0,3]\times[0,1]\), the locations of the troubled-cells at the final time and the blow-up region around the double Mach stems. The new hybrid HWENO scheme has better density resolutions than the  hybrid HWENO scheme, in addition, the hybrid HWENO scheme needs to use smaller CFL number taken as 0.45, but the CFL number for the new hybrid HWENO scheme is 0.6, moreover, the new hybrid HWENO scheme uses less candidate stencils but has higher order numerical accuracy.
\smallskip

\begin{figure}[!ht]
\centerline{\psfig{file=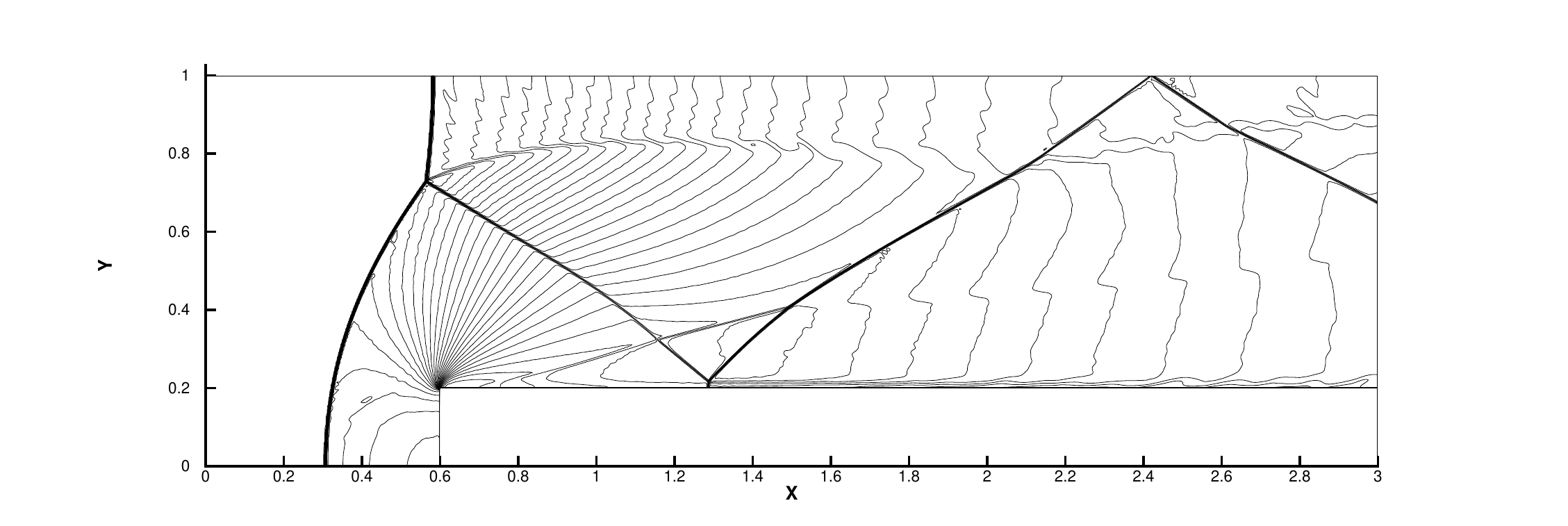,width=3.5 in} \psfig{file=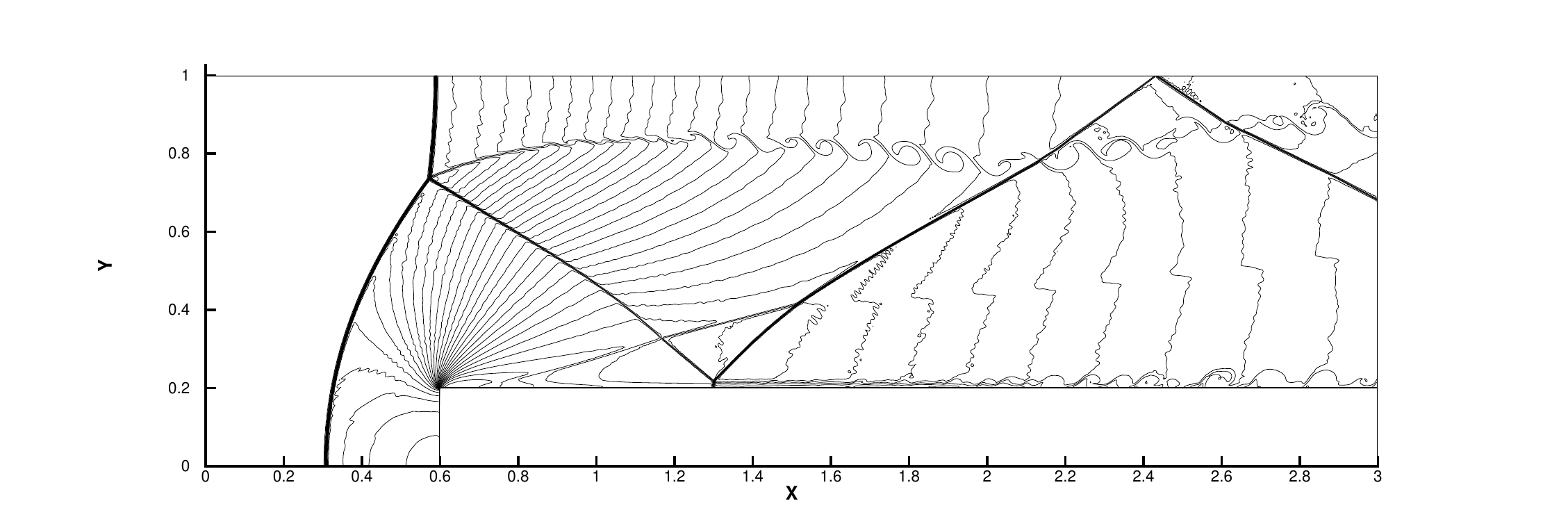,width=3.5 in} }
\centerline{\psfig{file=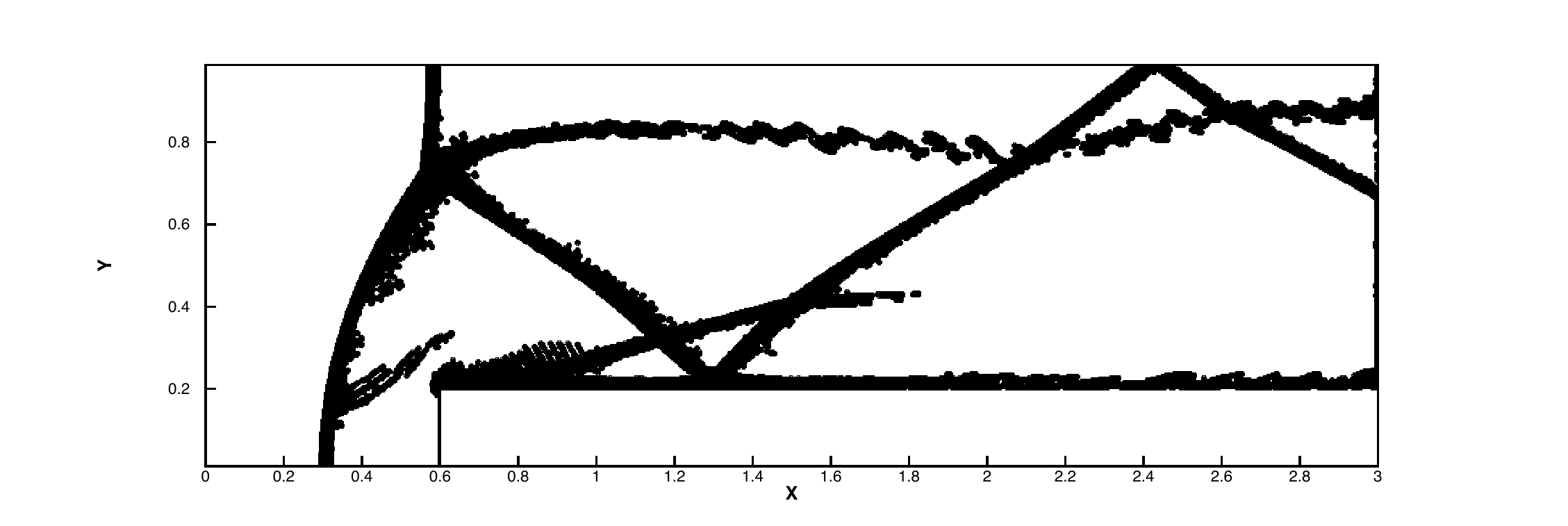,width=3.5 in} \psfig{file=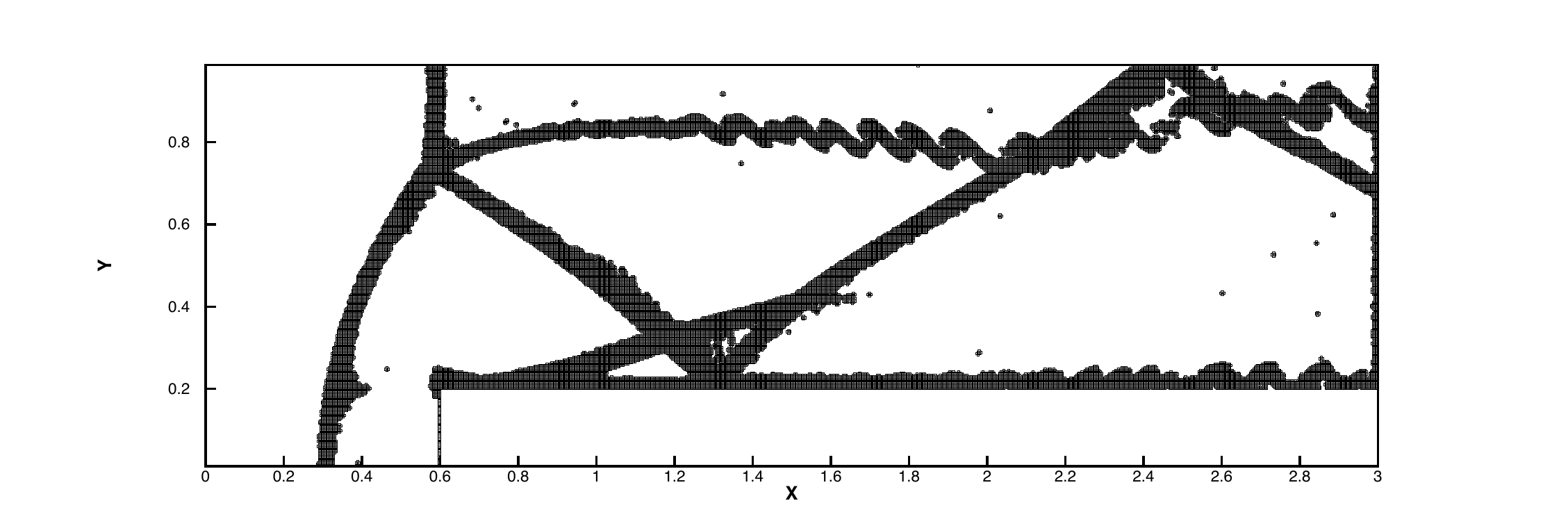,width=3.5 in} }
 \caption{ Forward step problem. T=4.
 From top to bottom: 30 equally spaced density contours from 0.32 to 6.15; the locations of the troubled-cells at the final time. The hybrid HWENO scheme (left); the new hybrid HWENO scheme (right). Uniform meshes with 960 $\times$ 320 cells.}
\label{stepfig}
\end{figure}

\noindent {\bf Example 3.11.} We finally solve the problem of a Mach 3 wind tunnel with a step \cite{Wooc} modeled by the two-dimensional Euler equations (\ref{euler2}). The wind tunnel is 1 length unit wide and 3 length units long. The step is 0.2 length units high and is located 0.6 length units from a right-going Mach 3 flow. Reflective boundary conditions are applied along the wall of the tunnel. In flow and out flow  boundary conditions are applied at the entrance and the exit, respectively. The computing time is up to $T=4$, then, we present the computed density and the locations of the troubled-cells at the final time in Figure \ref{stepfig}. We notice that the new hybrid HWENO scheme has high resolutions than the  hybrid HWENO scheme, and it also  has bigger CFL number, less candidate stencils, higher order numerical accuracy and simpler HWENO methodology. Similarly, only a small part of cells are identified as troubled-cells, and it means that most regions directly use linear approximation, which can increase the efficiency obviously.
\smallskip

\section{Concluding remarks}
\label{sec4}
\setcounter{equation}{0}
\setcounter{figure}{0}
\setcounter{table}{0}

In this paper, a  new fifth-order hybrid  finite volume Hermite weighted essentially non-oscillatory (HWENO) scheme with artificial linear weights is designed for solving  hyperbolic conservation laws. Compared with the hybrid HWENO scheme \cite{ZCQHH}, we employ a nonlinear convex combination of a high degree polynomial with several low degree polynomials in the new HWENO reconstruction, and the associated linear weights can be any artificial positive numbers (their sum is one), which would have the advantages of its simplicity and easy extension to multi-dimension. Meanwhile, different choice of the linear weights would not affect the numerical accuracy, and it gets less numerical errors than the original  HWENO methodology. In addition, the new hybrid HWENO scheme has higher order numerical accuracy in two dimension. Moreover, the scheme still keeps the non-oscillations as we apply the limiter methodology for the first order moments in the  troubled-cells and use new HWENO reconstruction on the interface. In the implementation, only a small part of cells are identified as troubled-cells, which means that most regions directly use linear approximation. In short, the new hybrid HWENO scheme has high resolution, efficiency, non-oscillation and robustness, simultaneously, and these numerical results also show  its good performances.



\begin{thebibliography}{99}

\bibitem{bgfb}
D. S. Balsara, S. Garain, V. Florinski and W. Boscheri, An efficient class of WENO schemes with adaptive order for unstructured meshes,  J. Comput. Phys., 404(2020), 109062.

\bibitem{bgs} D.S. Balsara, S. Garain and C.-W. Shu, An efficient class of WENO schemes with adaptive order, J. Comput. Phys., 326 (2016), 780-804.
\bibitem{DG2} B. Cockburn and C.-W. Shu, TVB Runge-Kutta local projection discontinuous Galerkin finite element method for conservation laws II: general framework, Mathematics of Computation, 52 (1989), 411-435.
\bibitem{cd1} B. Costa and W.S. Don,  Multi-domain hybrid spectral-WENO methods for hyperbolic conservation laws, J. Comput. Phys., 224 (2007), 970-991.
\bibitem{cd2}B. Costa and W.S. Don, High order Hybrid Central-WENO finite difference scheme for conservation laws, J. Comput. Appl. Math, 204 (2007), 209-218.

\bibitem{ccd} M. Castro, B. Costa and W.S. Don, High order weighted essentially non-oscillatory WENO-Z schemes for hyperbolic conservation laws, J. Comput. Phys., 230 (2011), 1766-1792.

\bibitem{CZQ} X. Cai. X. Zhang and J. Qiu, Positivity-preserving high order finite volume HWENO schemes for compressible Euler equations, J. Sci. Comput., 68 (2016), 464-483.

\bibitem{DBTM}
M. Dumbser, D.S. Balsara, E.F. Toro and C.D. Munz, A unified framework for the construction of one-step finite volume and discontinuous Galerkin schemes on unstructured meshes, J. Comput. Phys., 227 (2008), 8209-8253.

\bibitem{DBSR} M. Dumbser, W. Boscheri, M. Semplice and G. Russo, Central weighted ENO schemes for hyperbolic conservation laws on fixed and moving unstructured meshes, SIAM J. Sci. Comput., 39 (2017), A2564-A2591.

\bibitem{ho} A. Harten and S. Osher, Uniformly high-order accurate non-oscillatory
schemes, IMRC Technical Summary Rept. 2823, Univ. of Wisconsin,
Madison, WI, May 1985.
\bibitem{h1} A. Harten, Preliminary results on the extension of ENO schemes to two-dimensional problems, in Proceedings, International Conference on Nonlinear
Hyperbolic Problems, Saint-Etienne, 1986, Lecture Notes in Mathematics, edited by C. Carasso {\em et al.} (Springer-Verlag, Berlin, 1987).

\bibitem{heoc} A. Harten, B. Engquist, S. Osher and S. Chakravarthy, Uniformly high order accurate essentially non-oscillatory schemes III, J. Comput. Phys., 71 (1987), 231-323.

\bibitem{hs} C. Hu and C.-W. Shu, Weighted essentially non-oscillatory schemes on triangular meshes, J. Comput. Phys., 150 (1999), 97-127.

\bibitem{Hdp}D.J. Hill and  D.I. Pullin, Hybrid tuned center-difference-WENO method for large eddy simulations
in the presence of strong shocks, J. Comput. Phys.,  194 (2004), 435-450.

\bibitem{js} G.-S. Jiang and C.-W. Shu, Efficient implementation of weighted ENO schemes, J. Comput. Phys., 126 (1996), 202-228.

\bibitem{LJJN}L. Krivodonova, J. Xin, J.-F. Remacle, N. Chevaugeon and J.E. Flaherty, Shock detection and limiting with discontinuous Galerkin methods for hyperbolic conservation laws, Applied Numerical Mathematics, 48  (2004), 323-338.


\bibitem{loc} X.D. Liu, S. Osher and T. Chan, Weighted essentially
non-oscillatory schemes, J. Comput. Phys., 115 (1994), 200-212.

\bibitem{lpr} D. Levy, G. Puppo and G. Russo, Central WENO schemes for hyperbolic systems
of conservation laws, M2AN. Math. Model. Numer. Anal., 33 (1999), 547-571.

\bibitem{lpr2} D. Levy, G. Puppo and G. Russo, Compact central WENO schemes for
multidimensional conservation laws, SIAM J. Sci. Comput., 22 (2) (2000), 656-672.

\bibitem{Glj} G. Li and J. Qiu, Hybrid weighted essentially non-oscillatory schemes with different indicators, J. Comput. Phys., 229 (2010),  8105-8129.

\bibitem{LLQ}H. Liu and J. Qiu, Finite Difference Hermite WENO schemes for conservation laws,  J. Sci. Comput., 63 (2015), 548-572.

\bibitem{SPir}S. Pirozzoli, Conservative hybrid compact-WENO schemes for shock-turbulence interaction, J. Comput. Phys., 178 (2002), 81-117.

\bibitem{QSHw1}J. Qiu and C.-W. Shu, Hermite WENO schemes and their application as limiters for Runge-Kutta discontinuous Galerkin method: one-dimensional case, J. Comput. Phys., 193 (2004), 115-135.

\bibitem{QSHw}J. Qiu and C.-W. Shu, Hermite WENO schemes and their application as limiters for Runge-Kutta discontinuous Galerkin method II: Two dimensional case,  Computers \& Fluids, 34 (2005), 642-663.

\bibitem{so1} C.-W. Shu and S. Osher, Efficient implementation of
essentially non-oscillatory shock capturing schemes, J. Comput. Phys., 77 (1988), 439-471.

\bibitem{s2}
C.-W. Shu, Essentially non-oscillatory and weighted essentially non-oscillatory schemes for hyperbolic
conservation laws, In: Quarteroni, A. (ed.) Advanced Numerical Approximation of Nonlinear Hyperbolic
Equations, Lecture Notes in Mathematics, CIME subseries, Springer, Berlin (1998).

\bibitem{SHSN}J. Shi, C. Hu and C.-W. Shu, A technique of treating negative weights in WENO schemes, J. Comput. Phys., 175 (2002), 108-127.

\bibitem{TQ} Z. Tao, F. Li and J. Qiu, High-order central Hermite WENO schemes on staggered meshes for hyperbolic conservation laws, J. Comput. Phys., 281 (2015), 148-176.

\bibitem{TLQ}Z. Tao, F. Li and J. Qiu, High-order central Hermite WENO schemes: dimension-by-dimension moment-based reconstructions, J. Comput. Phys., 318 (2016), 222-251.

\bibitem{Wooc} P. Woodward and P. Colella, The numerical simulation of two-dimensional fluid flow with strong shocks. J. Comput. Phys., 54 (1984), 115-173.

\bibitem{ZQHW} J. Zhu and J. Qiu, A Class of Forth order Finite Volume Hermite Weighted Essentially Non-oscillatory Schemes,   Science in China, Series A--Mathematics,  51 (2008), 1549-1560.

\bibitem{ZQd} J. Zhu and J. Qiu, A new fifth order finite difference WENO scheme for solving hyperbolic conservation laws, J. Comput. Phys., 318 (2016), 110-121.

\bibitem{ZA}
 Y. H. Zahran and A. H. Abdalla,  Seventh order Hermite WENO scheme for hyperbolic conservation laws, Computers \& Fluids, 131 (2016), 66-80.

\bibitem{WZQ} J. Zhu and J. Qiu, A new fifth order finite difference WENO scheme for solving hyperbolic conservation laws, J. Comput. Phys., 318 (2016), 110-121.

\bibitem{vozq}J. Zhu and J. Qiu, A new type of finite volume WENO schemes for hyperbolic conservation laws, J. Sci. Comput., 73 (2017), 1338-1359.

\bibitem{tezq}
J. Zhu and J. Qiu, A new third order finite volume weighted essentially non-oscillatory scheme on tetrahedral meshes,   J. Comput. Phys., 349 (2017), 220-232.

\bibitem{trzq}
J. Zhu and J. Qiu, A simple finite volume weighted essentially non-oscillatory schemes on triangular meshes, SIAM J. Sci. Comput., 40 (2018), A903-A928.

\bibitem{ZZCQ}
Z. Zhao, J. Zhu, Y. Chen and J. Qiu,  A new hybrid WENO scheme for hyperbolic conservation laws,  Computers \& Fluids, 179 (2019), 422-436.

\bibitem{ZCQHH}
Z. Zhao, Y. Chen and J. Qiu,  A  hybrid Hermite WENO scheme for hyperbolic conservation laws,  J. Comput. Phys., 405 (2020), 109175.

\end{thebibliography}
\end{document}